\documentclass{gtart}

\def\ifplaintex{\expandafter\ifx\csname documentclass\endcsname\relax}

\ifplaintex 
\else
\fi

\expandafter\ifx\csname beginpicture\endcsname\relax
\expandafter\ifx\csname documentclass\endcsname\relax
\input pictex \else
\input prepictex \input pictex \input postpictex \fi\fi

\def\gt{{\mathsurround=0pt\it $\cal G\mskip-2mu$eometry \&\ 
$\cal T\!\!$opology}}        

\def\gtp{{\mathsurround=0pt\it $\cal G\mskip-2mu$eometry \&\ 
$\cal T\!\!$opology $\cal P\!$ublications}}  

\def\lognumber#1{\def\thelognumber{#1}}
\def\volumenumber#1{\def\thevolumenumber{#1}}
\def\papernumber#1{\def\thepapernumber{#1}}
\def\volumeyear#1{\def\thevolumeyear{#1}}

\def\pagenumbers#1#2{\def\startpage{#1}\def\finishpage{#2}}
\def\published#1{\def\publishdate{#1}}
\def\proposed#1{\def\theproposer{#1}}
\def\seconded#1{\def\theseconders{#1}}
\def\received#1{\def\receiveddate{#1}}
\def\revised#1{\def\reviseddate{#1}}
\def\accepted#1{\def\accepteddate{#1}}

\def\coverauthors#1{\def\thecoverauthors{#1}}
\def\asciiauthors#1{\def\theasciiauthors{#1}}
\def\asciiaddress#1{\def\theasciiaddress{#1}}
\def\asciiemail#1{\def\theasciiemail{#1}}

\let\\\par\let\thelognumber\relax
\let\thevolumenumber\relax\let\thepapernumber\relax
\let\thevolumeyear\relax\let\thesamplenumber\relax\let\startpage\relax
\let\finishpage\relax\let\publishdate\relax\let\receiveddate\relax
\let\reviseddate\relax\let\accepteddate\relax\let\theasciititle\relax
\let\theasciiauthors\relax\let\theasciiaddress\relax
\let\theasciiabstract\relax
\let\theasciiemail\relax\let\theshortauthors\relax\let\theshorttitle\relax
\let\thecoverauthors\relax

\long\def\maketitlep{   

\count0=\startpage

\gt\hfill      
\beginpicture
\setcoordinatesystem units <0.33truein, 0.33truein> point at 2.2 0.9
\setplotsymbol ({$\cal G$})
\plotsymbolspacing=9truept
\circulararc 315 degrees from 0 1 center at 0 0
\setplotsymbol ({$\cal T$})
\circulararc 315 degrees from 1 -1 center at 1 0
\endpicture
\break
{\small\ifx\thesamplenumber\relax 
Volume \else Sample
\fi\thevolumenumber\ (\thevolumeyear)
\startpage--\finishpage\nl
Published: \publishdate}
\vglue 0.5truein plus 0.4fil minus 0.1truein

{\parskip=0pt\leftskip 0pt plus 1fil\def\\{\par\smallskip}{\ifplaintex\large
\else\Large\fi\bf\thetitle}\par\medskip}   

\vglue 0pt plus 0.1fil 

{\parskip=0pt\leftskip 0pt plus 1fil\def\\{\par}{\sc\theauthors}
\par\medskip}

\vglue 0pt plus 0.1fil 

{\small\parskip=0pt\let\newline\\
{\leftskip 0pt plus 1fil\def\\{\par}{\sl\theaddress}\par}
\expandafter\ifx\theemail\relax    
\relax\else\vglue 5pt plus 0.02fil minus 2pt\def\\{\stdspace{\rm 
and}\stdspace} 
\cl{Email:\stdspace\tt\theemail}\fi
\ifx\theurl\relax                  
\relax\else\vglue 5pt plus 0.02fil minus 2pt\def\\{\stdspace{\rm 
and}\stdspace}
\cl{URL:\stdspace\tt\theurl}\fi\par}

\vglue 7pt plus 0.3fil minus 3pt

{\bf Abstract}
\vglue 5pt plus 0.1fil minus 2pt

\theabstract

\vglue 7pt plus 0.3fil minus 3pt

{\bf AMS Classification numbers}\quad Primary:\quad \theprimaryclass

Secondary:\quad \thesecondaryclass

\vglue 5pt plus 0.3fil minus 2pt

{\bf Keywords}\quad \thekeywords

\vglue 10pt plus 0.5fil minus 5pt

{\small  Proposed: \theproposer\hfill Received: \receiveddate\nl
Seconded: \theseconders\hfill 
\ifx\reviseddate\relax                         
Accepted: \accepteddate                        
\else
Revised: \reviseddate                          
\fi}
\eject
}       

\let\maketitlepage\maketitlep
\let\maketitle\maketitlepage

\font\phead=cmsl9 scaled 950
\font=cmsl9 scaled 1050
\font\pnum=cmbx10 scaled 913
\font\lnum=cmbx10 
\font\pfoot=cmsl9 scaled 950
\font=cmsl9 scaled 1050
\ifplaintex
\headline{\vbox to 0pt{\vskip -4.5mm\line{\small\phead\ifnum
\count0=\startpage ISSN 1364-0380 (on line)
1465-3060 (printed) \hfill {\pnum\folio}\else\ifodd\count0\def\\{ }%
\ifx\theshorttitle\relax\thetitle\else\theshorttitle\fi\hfill{\pnum\folio}
\else\def\\{ and }{\pnum\folio}\hfill\ifx\theshortauthors\relax\theauthors
\else\theshortauthors\fi\fi\fi}\vss}}
\footline{\vbox to 0pt{\vglue 0mm\line{\small\pfoot\ifnum\count0=\startpage
\copyright\ \gtp\hfill\else
\gt, Volume \thevolumenumber\ (\thevolumeyear)\hfill\fi}\vss
}}
\else
\makeatletter
\def\@oddhead{{\small\lhead\ifnum\count0=\startpage ISSN 1364-0380 (on line)
1465-3060 (printed) \hfill {\lnum\number\count0}\else\ifodd\count0
\def\\{ }\ifx\theshorttitle\relax \thetitle \else\theshorttitle\fi\hfill
{\lnum\number\count0}\else\def\\{ and }{\lnum\number\count0}
\hfill\ifx\theshortauthors\relax 
\theauthors\else\theshortauthors\fi\fi\fi}}\def\@evenhead{\@oddhead}
\def\@oddfoot{\small\lfoot\ifnum\count0=\startpage\copyright\ \gtp\hfill\else
\gt, Volume \thevolumenumber\ (\thevolumeyear)\hfill\fi}
\def\@evenfoot{\@oddfoot}
\makeatother
\fi

\newwrite\gtoutfile
\long\gdef\makeheadfile{  
{\def\\{, }\def\s{ }
\immediate\openout\gtoutfile head.xxx
\immediate\write\gtoutfile{To: math@arxiv.org}
\immediate\write\gtoutfile{Subject: put or rep NNNNN:pppp}
\immediate\write\gtoutfile{--text follows this line--}
\immediate\write\gtoutfile{Proxy-for: \ifx\theasciiauthors\relax
\theauthors\else\theasciiauthors\fi\s<\ifx\theasciiemail\relax\theemail\else\theasciiemail\fi>}
\immediate\write\gtoutfile{\noexpand\\}
\immediate\write\gtoutfile{Authors: \ifx\theasciiauthors\relax
\theauthors\else\theasciiauthors\fi}
{\def\\{ }\immediate\write\gtoutfile{Title: \ifx\theasciititle\relax
\thetitle\else\theasciititle\fi}}
\immediate\write\gtoutfile{Subj-class: GT or SG or MG etc}
\immediate\write\gtoutfile{MSC-class: \theprimaryclass\ifx\thesecondaryclass\relax\else, \thesecondaryclass\fi}
\immediate\write\gtoutfile{Journal-ref: Geom. Topol. \thevolumenumber
(\thevolumeyear) \startpage-\finishpage}
\immediate\write\gtoutfile{Comments: Published by Geometry and Topology at}
\immediate\write\gtoutfile{\s\s http://www.maths.warwick.ac.uk/gt/GTVol\thevolumenumber/paper\thepapernumber.abs.html}
\immediate\write\gtoutfile{\noexpand\\}
\immediate\write\gtoutfile{}
\ifx\theasciiabstract\relax
\immediate\write\gtoutfile{\theabstract}\else
\immediate\write\gtoutfile{\theasciiabstract}\fi
\immediate\write\gtoutfile{}
\immediate\write\gtoutfile{\noexpand\\}
\immediate\write\gtoutfile{}
\immediate\closeout\gtoutfile}}  

\def\maketitlepage{\maketitlep\makeheadfile}
\let\maketitle\maketitlepage

\lognumber{247}
\volumenumber{7}\papernumber{5}\volumeyear{2003}
\pagenumbers{185}{224}
\received{15 April 2002}
\revised{28 January 2003}
\accepted{9 March 2003}
\published{24 March 2003}
\proposed{Robion Kirby}
\seconded{Tomasz Mrowka, Yasha Eliashberg}

\usepackage{amsmath,amsfonts,amscd,flafter,rlepsf}

\newtheorem{theorem}{Theorem}[section]
\newtheorem{prop}[theorem]{Proposition}
\newtheorem{cor}[theorem]{Corollary}
\newtheorem{conj}[theorem]{Conjecture}
\newtheorem{lemma}[theorem]{Lemma}

\theoremstyle{remark}
\newtheorem{defn}[theorem]{Definition}

\newtheorem{remark}[theorem]{Remark}

\newcommand\Rk{\mathrm{rk}}

\newcommand{\HF}{HF}

\newcommand{\Q}{\mathbb{Q}}

\newcommand{\Z}{\mathbb{Z}}

\newcommand{\Zmod}[1]{\Z/{#1}\Z}

\newcommand{\Ker}{\mathrm{Ker}}

\newcommand{\cm}{\cdot}

\newcommand{\ModSWfour}{\mathcal{M}}
\newcommand{\ModFlow}{\ModSWfour}

\newcommand{\SpinC}{{\mathrm{Spin}}^c}

\newcommand\Wedge{\Lambda}

\newcommand\Hom{\mathrm{Hom}}

\newcommand\abuts\Rightarrow
\newcommand\Sym{\mathrm{Sym}}

\newcommand\mCP{{\overline{\mathbb{CP}}}^2}

\newcommand\HFpRed{\HFp_{\red}}

\newcommand{\ev}{\mathrm{ev}}
\newcommand{\odd}{\mathrm{odd}}

\newcommand\Filt{\mathcal F}
\newcommand\HFinfty{\HFinf}

\newcommand\ModSphere{\ModFlow\left({\mathbb S}\longrightarrow 
\Sym^{g-1}(\Sigma_{1})\times \Sym^2(\Sigma_{2})\right)}
\newcommand\ModSpheres\ModSphere

\newcommand{\red}{\mathrm{red}}

\newcommand\HFp{\HFb}

\newcommand\HFm{\HF^-}

\newcommand\HFinf{HF^\infty}

\newcommand\HFa{\widehat{HF}}
\newcommand\HFb{HF^+}

\newcommand\UnparModSp{\widehat \ModSp}
\newcommand\UnparModFlow\UnparModSp
\newcommand\Mod\ModSp

\newcommand\PD{\mathrm{PD}}

\newcommand{\spinc}{\mathfrak s}

\newcommand{\spinct}{\mathfrak t}

\newcommand\ModMaps{\mathcal M}
\newcommand\ModSp\ModMaps

\newcommand\uHF{\underline{\HF}}

\newcommand\uHFinf{\uHF^\infty}

\newcommand\Fp[1]{F^{+}_{#1}}

\newcommand\Finf[1]{F^{\infty}_{#1}}

\newcommand\Field{\mathbb F}

\newcommand\Dual{\mathcal D}
\newcommand\Duality\Dual

\newcommand\Knot{\mathbb K}

\newcommand\Vertices{\mathrm{Vert}}
\newcommand\Char{\mathrm{Char}}
\newcommand\Combp{{\mathbb H}^+}
\newcommand\DCombp{{\mathbb K}^+}
\newcommand\Comba{\widehat{\mathbb H}}
\newcommand\NatTransp{T^+}
\newcommand\InjMod[1]{{\mathcal T}^+_{#1}}
\newcommand\NatTransa{\widehat T}
\newcommand\MapOnep{A^+}
\newcommand\MapTwop{B^+}
\newcommand\MapThreep{C^+}
\newcommand\MapOnea{\widehat A}
\newcommand\MapTwoa{\widehat B}
\newcommand\MapOneCombp{{\mathbb A}^+}
\newcommand\MapTwoCombp{{\mathbb B}^+}
\newcommand\MapOneComba{\widehat {\mathbb A}}
\newcommand\MapTwoComba{\widehat {\mathbb B}}

\newcommand\fp[1]{f^+_{#1}}
\newcommand\Short{\mathcal S}

\begin{document}

\title{On the Floer homology of plumbed three-manifolds}

\author{Peter Ozsv\'ath\\Zolt{\'a}n Szab{\'o}}
\coverauthors{Peter Ozsv\noexpand\'ath\\Zolt{\noexpand\'a}n Szab{\noexpand\'o}}
\asciiauthors{Peter Ozsvath, Zoltan Szabo} 
\address{Department of
Mathematics, Columbia University\\New York 10027, USA}
\email{petero@math.columbia.edu}
\secondaddress{Department of
Mathematics, Princeton University\\New Jersey 08540, USA}
\secondemail{szabo@math.princeton.edu}

\asciiaddress{Department of
Mathematics, Columbia University\\New York 10027, USA\\and\\Department of
Mathematics, Princeton University\\New Jersey 08540, USA}
\asciiemail{petero@math.columbia.edu, szabo@math.princeton.edu}

\begin{abstract}  
We calculate the Heegaard Floer homologies for three-manifolds
obtained by plumbings of spheres specified by certain graphs. Our
class of graphs is sufficiently large to describe, for example, all
Seifert fibered rational homology spheres. These calculations can be
used to determine also these groups for other three-manifolds,
including the product of a circle with a genus two surface.
\end{abstract}


\maketitlepage

\section{Introduction}

In~\cite{HolDisk}, we defined Heegaard Floer homology invariants for
closed, oriented three-manifolds. In~\cite{HolDiskFour}, we defined
invariants for cobordisms between three-manifolds, and consequently
also for smooth, closed four-manifolds.  The resulting package has
many properties of a topological quantum field theory, and moreover it
is closely related to its gauge-theoretic counterparts,
Donaldson-Floer (see~\cite{DonaldsonFloer}) and Seiberg-Witten theory
(see~\cite{Witten}, \cite{Morgan}, \cite{KMthom}).

In particular, the three-manifold invariants are a fundamental
stepping-stone in the definition and computation of the four-manifold
invariants.  Moreover, many four-dimensional aspects of the
three-manifold $Y$ are reflected in its Heegaard Floer homology,
including obstructions to embedding the three-manifold in a symplectic
four-manifold (cf.\ \cite{HolDiskSymp}) and also restrictions on the
intersection forms of smooth four-manifolds which bound $Y$
(see~\cite{AbsGraded}, compare~\cite{Froyshov}).

Whereas ingredients in the Heegaard Floer homology are more
combinatorial in flavor than the corresponding gauge theory
ingredients, the definition still involves a fundamentally analytical
object: holomorphic disks in the symmetric product of a Riemann
surface. Our aim here is to give a combinatorial formulation of these
groups for a class of three-manifolds which are obtained by certain 
plumbing
diagrams. Indeed, this class is large enough to describe, for example,
all Seifert fibered rational homology spheres. The answer we describe
can be read off from the plumbing tree. Moreover, this answer
contains, as a by-product, all of the relative invariants of the
four-manifold obtained from the plumbing description. It is
interesting to compare these calculations with their corresponding
analogues in instanton Floer homology (see for
example~\cite{FintushelStern}, \cite{KirkKlassen}) and Seiberg-Witten
theory (see for example~\cite{MOY}). Note that, for many of the
three-manifolds studied in this paper, the corresponding instanton
Floer homology and Seiberg-Witten theory remain elusive.

Applications of these calculations include, as we have mentioned,
non-embedding theorems for certain of these three-manifolds in
symplectic four-manifolds, cf.\ \cite{HolDiskSymp}. As another
application, we show that the four-manifold obtained from the plumbing
description has, in some sense, a maximally exotic intersection form,
as measured by the lengths of characteristic vectors,
cf.\ Corollary~\ref{cor:IntersectionForm} below.  The computations in
this paper also play a major role in~\cite{SeifSurg}, where
we give constraints on knots in the three-sphere which admit Seifert
fibered surgeries. A final application described here gives
calculations of the Heegaard Floer homology for some other
three-manifolds, including the product of a circle with a surface of
genus two. This latter calculation is used to shed some light on the
structure of Heegaard Floer homology for more complicated
three-manifolds.

With this motivation in hand, we describe the family of
three-manifolds studied in this paper; but first, we give some
preliminaries.

We call a {\em weighted graph} $G$ a graph equipped with an
integer-valued function $m$ on its vertices.  A weighted graph
gives rise to a four-manifold with boundary $X(G)$ which is obtained
by plumbing together a collection of disk bundles over the two-sphere
(indexed by vertices of $G$), so that the Euler number of the sphere
bundle corresponding to the vertex $v$ is given by its multiplicity
$m(v)$.  The sphere belonging
to $v$ is plumbed to the sphere belonging to $w$ precisely when the
two are connected by an edge. Let $Y(G)$ be the oriented
three-manifold which is the boundary of $X(G)$. 

For $X=X(G)$, the group $H_2(X;\Z)$ is the lattice freely spanned by
the vertices of $G$, and the intersection form on $H_2(X;\Z)$ is given
by the graph as follows. For a vertex $v$ of $G$, let $[v]\in
H_2(X;\Z)$ denote the corresponding homology class.  Then, for each
vertex $[v]\cm [v]=m(v)$, and for each pair of distinct vertices $v$ and $w$,
$[v]\cm [w]$ is one if $v$ and $w$ are connected by an
edge, and zero otherwise.

\begin{defn}
\label{def:OneBadPoint}
A weighted graph is said to be a {\em negative-definite graph} if:
\begin{itemize}
\item $G$ is a disjoint union of trees
\item the intersection form associated to $G$ is negative definite.
\end{itemize}
The {\em degree} of a vertex $v\in{\mathrm{Vert}}(G)$, denoted $d(v)$,
is the number of edges which contain $v$.
A vertex $v\in {\mathrm{Vert}}(G)$ is said to be a {\em bad vertex} of the weighted graph if
$$m(v)>-d(v).$$
\end{defn}

In this paper, we will be primarily concerned with negative-definite
graphs with at most one bad vertex.

Note that any Seifert fibered rational homology sphere (with at least
one orientation) can be realized from a negative-definite graph which
is star-like (i.e.\ is a connected graph with at most one vertex with
degree $>2$), so that if $v$ is a vertex with degree $d(v)\leq 2$,
then $m(v)\leq -2$ (see for example~\cite{GompfStipsicz}). In particular,
this is a negative-definite graph with at most one bad vertex.

Our goal here is to give an algebraic description of the Heegaard
Floer homology groups $\HFp(-Y(G))$. Recall that the Heegaard Floer
homology groups come in a package, $\HFp$, $\HFm$, $\HFinf$ and $\HFa$
which are all closely related. However, for a rational homology
three-sphere, all of the information can be extracted from
$\HFp$. Recall that this group is in general a module over the ring
$\Z[U]$, where $U$ lowers degree by two, and every element in $\HFp$
is annihilated by a sufficiently large power of $U$. 

As a
starting point, let
$\InjMod{0}$ denote the graded $\Z[U]$-module which is the quotient of 
$\Z[U,U^{-1}]$ by the submodule $U\cm \Z[U]$. 
This module is graded so that the element 
$U^{-d}$ (for $d\geq 0$) is supported in degree $2d$.
Recall that
$$\HFp(S^3)\cong \InjMod{0}.$$
Let $\Char(G)$ denote the set of characteristic vectors for the
intersection form. Let 
$$\Combp(G)\subset \Hom(\Char(G),\InjMod{0})$$
denote the set of functions with finite support and 
which satisfy the following ``adjunction relations''
for all characteristic vectors $K$ and vertices $v$. 
Let
$$2n=\langle K,v \rangle + v\cm v.$$
If $n\geq 0$, then we require that
\begin{equation}
\label{eq:AdjRel}
U^n\cm \phi(K+2\PD[v]) = \phi(K),
\end{equation}
while if $n\leq 0$, then
\begin{equation}
\label{eq:AdjRel2}
\phi(K+2\PD[v]) = U^{-n}\cm \phi(K).
\end{equation}
We can decompose $\Combp(G)$ according to $\SpinC$ structures over $Y$.
Note first that the first Chern class gives an identification of the set of
$\SpinC$ structures over $X=X(G)$ with the set of characteristic vectors
$\Char(G)$. Observe that the image of $H^2(X,\partial X;\Z)$ in $H^2(W;\Z)$
is spanned by the Poincar\'e duals of the spheres corresponding to the vertices.
Using the restriction to boundary, it is easy to see that the set of $\SpinC$ structures
over $Y$ is identified with the set of $2 H^2(X,\partial X;\Z)$-orbits 
in $\Char(G)$.

Fix a $\SpinC$ structure $\spinct$ over $Y$. Let $\Char_\spinct(G)$
denote the set of characteristic vectors for $X$ which are first
Chern classes of $\SpinC$ structures $\spinc$ whose restriction to the
boundary is $\spinct$. Similarly, we let $$\Combp(G,\spinct)\subset
\Combp(G)$$ be the subset of maps which are supported on the subset of
characteristic vectors $\Char_\spinct(G)\subset \Char(G)$.  We have a
direct sum splitting: $$\Combp(G)\cong \bigoplus_{\spinct\in\SpinC(Y)}
\Combp(G,\spinct).$$
We can also introduce a grading on $\Combp(G)$ as follows.
We say that an element $\phi\in \Combp(G)$ is homogeneous of degree $d$ if 
for each characteristic vector $K$ with $\phi(K)\neq 0$, $\phi(K)\in \InjMod{0}$ is 
a homogeneous element with:
\begin{equation}
\label{eq:DefOfDegree}
\deg(\phi(K))-\left(\frac{K^2+|G|}{4}\right)=d.
\end{equation}
Our main result is the following identification of $\HFp(-Y(G))$ in
terms of combinatorics of the plumbing diagram:

\begin{theorem}
\label{intro:SomePlumbings}
Let $G$ be a negative-definite weighted graph with at most one bad vertex, in the sense
of Definition~\ref{def:OneBadPoint}. Then, for each $\SpinC$ structure
$\spinct$ over $-Y(G)$, there is an isomorphism of graded
$\Z[U]$ modules,
$$\HFp(-Y(G),\spinct)\cong \Combp(G,\spinct).$$
\end{theorem}

\begin{remark}
It is a straightforward matter to determine $\HFp(Y(G),\spinct)$
from $\HFp(-Y(G),\spinct)$, cf.\ Section 2 of [12].
\end{remark}

In the statement of the above theorem, the grading on
$\HFp(-Y(G),\spinct)$ is the absolute $\Q$-grading defined
in~\cite{HolDiskFour} and studied in~\cite{AbsGraded}.  Recall that when
$-Y(G)$ is an integral homology sphere, this absolute grading takes
values in $\Z$.

As a qualitative remark, it is perhaps worth pointing out the
following corollary (compare~\cite{FintushelStern}). To state it,
recall that there is an absolute $\Zmod{2}$-grading on $\HFp(Y,\spinct)$
which, for rational homology three-spheres, is determined by the
following criterion.  A homogeneous element $\xi$ is even with respect
to this grading if there is a non-zero homogeneous element of
$\xi_0\in\HFinf(Y,\spinct)$ with the property that
$$\deg(\xi)-\deg(\xi_0)\equiv 0\pmod{2}.$$ When $Y$ is an integral
homology sphere, this notion coincides with the parity of the
($\Z$-)grading of $\xi$.

\begin{cor}
\label{cor:EvenDegrees}
If $G$ is a negative-definite graph with at most one 
bad vertex, then all
elements of $\HFp(-Y(G),\spinct)$ have even $\Zmod{2}$ grading.
\end{cor}

\begin{proof}
This follows immediately from Theorem~\ref{intro:SomePlumbings}
and the definition of the absolute gradings: $\phi(K)\in\InjMod{0}$,
and the latter module is supported only in even degrees.
\end{proof}

This underscores the importance of the hypothesis on the graph. For
example, if $Y$ is the Brieskorn homology sphere $\Sigma(2,3,7)$
(which can be thought of as $(-1)$-surgery on the right-handed trefoil
knot) then it follows easily from the K\"unneth formula for connected
sums (Theorem~\ref{HolDiskTwo:thm:ConnSumHFm} of~\cite{HolDiskTwo})
that $\HFp(-(Y\#Y))$ has elements of both parities. On the other hand,
$Y\#Y$ admits a plumbing description as a negative-definite
disconnected graph with two bad points. For an example belonging to a
connected graph, one can take $-1$ surgery on the connected sum of two
right-handed trefoil knots in $S^3$, see
Proposition~\ref{prop:DoubleTrefoil} and Remark~\ref{rmk:OddGens}
below.  The methods for obtaining Theorem~\ref{intro:SomePlumbings}
do, however, give information on the Floer homology groups of the
three-manifolds obtained from these plumbing diagrams as well, see
Theorem~\ref{thm:PushedFurther} below.

Theorem~\ref{thm:SomePlumbings} also has the following corollary. For
the purpose of this corollary, recall that in~\cite{AbsGraded}, we
defined an invariant $d(Y,\spinct)$ associated to an oriented,
rational homology three-sphere $Y$ equipped with a $\SpinC$ structure
$\spinct$. This invariant takes values in $\Q$, The importance of
$d(Y,\spinct)$ is shown by the fact that it gives a bound on the
exoticness of the intersection form for any smooth, definite
four-manifold which bounds $Y$. Specifically, if $Y$ is a rational
homology three-sphere equipped with a $\SpinC$ structure $\spinct$,
then if $W$ is an oriented four-manifold with negative-definite
intersection form, and $\spinc$ is any $\SpinC$ structure over $W$
whose restriction to $Y$ is $\spinct$, then
Theorem~\ref{AbsGraded:thm:IntFormQSphere} of~\cite{AbsGraded}
establishes the inequality
\begin{equation}
\label{ineq:BoundShadow}
c_1(\spinc)^2+\Rk H^2(X;\Z)\leq 4d(Y,\spinct).
\end{equation}
Compare also the gauge-theoretic version of Fr{\o}yshov,
~\cite{Froyshov} and~\cite{FroyshovII}. (For the relationship between
diagonalizability of definite, unimodular forms $Q$ and the maximal
value, over all characteristic vectors $K$ for $Q$, of the quantity
$K^2+\Rk$, see~\cite{Elkies}.)  
We have the following consequence of Theorem~\ref{thm:SomePlumbings}
(which, in the case where $G$ has two bad points, follows from
Theorem~\ref{thm:PushedFurther}):

\begin{cor}
\label{cor:CalcD}
Let $G$ be a negative-definite graph with at most two bad points, and fix a $\SpinC$ structure
$\spinct$ over $Y$. Then, 
\begin{equation}
\label{eq:IdentifyDs}
d(Y(G),\spinct)=
\max_{\{K\in\Char_\spinct(G)\}}
\frac{K^2+|G|}{4}.
\end{equation}
\end{cor}

The above result gives a practical calculation of $d(Y,\spinct)$: for
a given $\spinct\in\SpinC(Y)$, it is easy to see that the maximum of
$\frac{K^2+|G|}{4}$ is always achieved among the finitely many
characteristic vectors $K\in\Char_\spinct(G)$ with $$|K\cm v| \leq
|m(v)|.$$ (A smaller set containing these minimal vectors is described
in Proposition~\ref{prop:KerU} below.)

Inequality~\eqref{ineq:BoundShadow},
combined with Corollary~\ref{cor:CalcD}, immediately gives the following:

\begin{cor}
\label{cor:IntersectionForm}
Let $G$ be a negative-definite graph with at most two bad points, and fix a $\SpinC$ structure
$\spinct$ over $Y$.
Then, for each smooth, compact, oriented four-manifold $X$ 
with negative intersection form which bounds $Y$, and for each
$\SpinC$ structure $\spinc\in\SpinC(X)$ with $\spinc|Y=\spinct$, we
have that $$c_1(\spinc)^2+\Rk (H^2(X;\Z))\leq
\max_{\{K\in\Char_\spinct(G)\}}
K^2+|G|.$$
\end{cor}

The above results are proved in Section~\ref{sec:Proof}. In
Section~\ref{sec:Examples} we give some sample calculations. In
Section~\ref{sec:GenusTwo}, we use these techniques as a
starting-point for another calculation: the calculation of
$\HFp(S^1\times\Sigma_2)$ (cf.\ Theorem~\ref{thm:SOneSigmaTwo} below).

We end the paper with some speculations based on this latter result.
Specifically, recall that we defined in~\cite{HolDisk} and
\cite{HolDiskTwo} a group $\HFinf$ which captures the behaviour
of $\HFp$ in all sufficiently large degrees. When the three-manifold
has $b_1(Y)<3$, $\HFinf$ is determined by $b_1(Y)$. It remains an
interesting question to determine $\HFinf$ for arbitrary
three-manifolds. We conclude this paper with a conjecture relating
$\HFinf(Y)$ with the cohomology ring of $Y$.

\medskip
\noindent{\bf{Acknowledgements}}\qua We would like to thank the 
referee for a careful reading of the manuscript and some very useful
comments.

PSO was supported by NSF grant number DMS 9971950 and a Sloan 
Research Fellowship; ZSz was supported by NSF grant number DMS 0107792
and a Packard Fellowship.

\section{Proof of Theorem~\ref{intro:SomePlumbings}}
\label{sec:Proof}

In the present section, we state a more precise version of
Theorem~\ref{intro:SomePlumbings}, and give a proof.  For the more
precise statement, we need the following notions.

First, to postpone a discussion of signs which might obscure matters,
we work over the field with two elements $\Field=\Zmod{2}$ for the
rest of the subsection, returning to a sign-refinement which allows us
to work over $\Z$ in Subsection~\ref{subsec:Signs}. Thus, unless it is
explicitly stated otherwise, all Floer homology groups in this subsection
are meant to be
taken with $\Field$ coefficients (which we suppress from the
notation).  In particular, with these conventions, $\InjMod{0}$ now
denotes the quotient of $\Field[U,U^{-1}]$ by the submodule $U\cm
\Field[U]$.

We define a map $$\NatTransp\colon \HFp(-Y(G))\longrightarrow
\Combp(G)$$ as follows. The plumbing diagram can be viewed as giving a
cobordism $W(G)$ from $S^3$ to the three-manifold $Y(G)$ (i.e.\  this
is the four-manifold obtained by deleting a ball from the
four-manifold $X(G)$ considered in the introduction) or, equivalently,
a cobordism from $-Y(G)$ to $S^3$. Now let $$\NatTransp(\xi)\colon
\Char(G)\longrightarrow \InjMod{0}$$ be the map given by
$$\NatTransp(\xi)(K)=\Fp{W(G),\spinc}(\xi)\in\HFp(S^3)=\InjMod{0},$$ where
$\spinc\in\SpinC(W(G))$ is the $\SpinC$ structure whose first Chern
class is $K$, and $\Fp{W(G),\spinc}$ denotes the four-dimensional
cobodism invariant defined in~\cite{HolDiskFour}.

\begin{theorem}
\label{thm:SomePlumbings}
Let $G$ be a negative-definite graph with at most one bad vertex. Then,
$T^+$ induces a grading-preserving isomorphism:
$$\Combp(G,\spinct)\cong\HFp(-Y(G),\spinct).$$
\end{theorem}

These techniques can be pushed further to obtain the following:

\begin{theorem}
\label{thm:PushedFurther}
Let $G$ be a negative-definite graph with at most two bad vertices.
Then, $T^+$ produces an isomorphism of graded $\Z[U]$-modules
$$\Combp(G,\spinct)\cong\HFp_{\ev}(-Y(G),\spinct),$$
where $\HFp_{\ev}$ denotes the part of $\HFp$ with even parity
(using the absolute $\Zmod{2}$ grading).
\end{theorem}

In practice, it is sometimes easier to think about $\Combp(G)$ from
the following dual point of view. We let $\DCombp(G)$ denote the
equivalence classes in $\Z^{\geq 0}\times
\Char(G)$ (where we write a pair $m$ and $K$ as $U^m\otimes K$) under 
the following equivalence relation. Let $v$ be a vertex and let
$$2n=\langle K,v \rangle + v\cm v.$$
If $n\geq 0$, then:
$$U^{n+m}\otimes (K+2\PD[v]) \sim U^m\otimes K,$$
while if $n\leq 0$, then
$$U^m\otimes (K+2\PD[v]) \sim U^{m-n}\otimes K.$$
Given a function $$\phi\colon \Char(G)\longrightarrow \InjMod{0},$$
there is an induced map $${\widetilde \phi}\colon \Z^{\geq 0}\times
\Char(G)\longrightarrow \InjMod{0}$$
defined by $${\widetilde \phi}(U^n\otimes K)=U^n\cm \phi(K).$$
Clearly, the set of finitely-supported functions $\phi\colon \Char(G)\longrightarrow
\InjMod{0}$ whose induced map ${\widetilde \phi}$ descends to
$\DCombp(G)$ is precisely $\Combp(G)$.

\begin{lemma}
\label{lemma:Dualize}
Let ${\mathcal B}_n$ denote the set of characteristic vectors 
$${\mathcal B}_n=\{K\in\Char(G) \big| \forall v\in G, |\langle K,v\rangle |\leq -m(v)+2n\}.$$
The quotient map induces a surjection from 
$$\bigcup_{i=0}^n U^i\otimes {\mathcal B}_{n-i}$$
onto the quotient space
$$\frac{{\DCombp}(G)}{\Z^{>n}\times\Char(G)}.$$ 
In turn, we have an identification
$$\big(\Ker U^{n+1}\subset \Combp(G;\Field)\big)
\cong \Hom\left(
\frac{{\DCombp}(G)}{\Z^{>n}\times\Char(G)},\Field\right)$$
(i.e.\ the right-hand-side consists of maps from $\DCombp(G)$ to $\Field$
which vanish on the equivalence classes which contain representatives of the
form $U^m\otimes K'$ with $m>n$)
and, indeed, 
$$\big(\Ker U^{n+1}\subset \Combp(G;\Z)\big)
\cong \Hom\left(
\frac{{\DCombp}(G)}{\Z^{>n}\times\Char(G)},\Z\right).$$
\end{lemma}

\begin{proof}
The surjectivity statement follows easily from the definition of the
equivalence relation in $\DCombp(G)$. 

The duality map is the one sending $$\phi\times (U^\ell\otimes K)
\mapsto \left(U^\ell \cm \phi(K)\right)_0$$ (i.e.\ taking the part in
$\InjMod{0}$ which lies in degree zero).  This obviously induces a map
$$ \Ker U^{n+1} \longrightarrow
\Hom\left(
\frac{{\DCombp}(G)}{\Z^{>n}\times\Char(G)},\Z\right).$$ This map is
injective, since if $\phi(K)\in \InjMod{0}$ is an element with
$$\left(U^\ell \cm \phi(K)\right)_0\equiv 0$$ for all $\ell>n$, then
clearly $U^{n+1}\cm \phi(K)=0$. To see that the map is surjective, observe that
if $$\tau\in
\Hom\left(
\frac{{\DCombp}(G)}{\Z^{>n}\times\Char(G)},\Z\right)$$ is an
arbitrary element, we can define a map $$\phi(K)=\sum_{\ell=0}^n
\tau(U^\ell\otimes K)\cm U^{-\ell} $$ whose degree zero part is $\tau$. 
Clearly, $\phi\in\Combp(G)$, and $U^{n+1}\cm \phi=0$.
\end{proof}

We now set up some properties of $\Combp(G,\spinct)$ with a view
towards proving Theorem~\ref{thm:SomePlumbings}.

\begin{prop}
\label{prop:SetUp}
The map $\NatTransp$ induces an $\Field[U]$-equivariant, degree-preserving
map from $\HFp(-Y(G),\spinct)$ to $\Hom(\Char_{\spinct}(G),\InjMod{0})$ whose image lies
in $$\Combp(G,\spinct)\subset \Hom(\Char_\spinct(G),\InjMod{0}).$$ 
\end{prop}

\begin{proof}
The map $T^+$ lands inside $\Hom(\Char_\spinct(G),\InjMod{0})$ with finite
support, according to general finiteness properties of the maps on
$\HFp$ induced by cobordisms
(cf.\ Theorem~\ref{HolDiskFour:thm:Finiteness}
of~\cite{HolDiskFour}). Alternatively, this finiteness follows from
the degree shift formula for maps induced by cobordisms, 
Theorem~\ref{HolDiskFour:thm:AbsGrade} of~\cite{HolDiskFour}, which
also shows that $T^+$ is degree-preserving. The fact that
$T^+$ lands in the subset of $\Hom(\Char_\spinct(G),\InjMod{0})$ satisfying
the adjunction relation (Equation~\eqref{eq:AdjRel}
or~\eqref{eq:AdjRel2} as appropriate) defining $\Combp(G)$ is proved in
Theorem~\ref{HolDiskSymp:thm:AdjunctionRelation} of~\cite{HolDiskSymp}
(where, actually, relations are established for oriented, embedded
surfaces of arbitrary genus).
\end{proof}

If $G$ is a weighted graph with a distinguished vertex
$v\in\Vertices(G)$, we let $G'(v)$ be a new graph formed by
introducing one new vertex $e$ labelled with weight $-1$, and
connected to only one other vertex, $v$.  Moreover, we let $G_{+1}(v)$
denote the weighted graph whose underlying graph agrees with $G$, but
whose weight at $v$ is increased by one (and the weight stays the same
for all other vertices). 
The two three-manifolds $Y(G'(v))$ and $Y(G_{+1}(v))$ are
clearly diffeomorphic; in fact, we have the following:

\begin{prop}
\label{prop:BlowDown}
Let $G'(v)$ be the graph obtained from $G$ as above. Then, there is a
grading-preserving isomorphism $$R\colon \Combp(G'(v))\rightarrow
\Combp(G_{+1}(v))$$ Moreover, this map is natural with $\NatTransp$ in
the sense that: $$\begin{CD}
\HFp(-Y(G'(v))) @>{\Theta}>> \HFp(-Y(G_{+1}(v))) \\
@V{T^+_{G'(v)}}VV @V{T^+_{G_{+1}(v)}}VV \\
\Combp(G'(v)) @>{R}>> \Combp(G_{+1}(v)),
\end{CD}$$
where map $\Theta$ is the isomorphism induced by the diffeomorphism of
$Y(G'(v))\cong Y(G_{+1}(v))$.
\end{prop}

\begin{proof}
We construct $R$ in two steps. As a first step, let
$G_{+1}(v)\cup f$ denote the disconnected graph consisting of the
disjoint union of $G_{+1}(v)$ and a single vertex $f$ with
multiplicity $-1$.  We have a map $$\Char(G_{+1}(v)\cup f)\cong
\Char(G'(v))$$ induced by a change of basis. It is easy to see that
this map induces an isomorphism $$\Combp(G_{+1}(v)\cup f)\cong
\Combp(G'(v)).$$
Next, we define a map
$$Q\colon \Combp(G_{+1}(v))\longrightarrow \Combp(G_{+1}(v)\cup f)$$
by the formula (where $m\geq 0$)
\begin{eqnarray*}
Q(\phi)(K,(2m+1))&=&U^{m(m+1)/2}\cm \phi(K), \\
Q(\phi)(K,-(2m+1))&=&U^{m(m+1)/2}\cm \phi(K).
\end{eqnarray*}
In the above notation, $(K,\ell)$ denotes the characteristic vector for
$G_{+1}(v)\cup f$ whose restriction to $G_{+1}(v)$ is $K$ and whose evaluation
on $f$ is $\ell$.
The map $Q$, too, is clearly an isomorphism. 
We define $R$ to be the composition of the above isomorphisms.

To check commutativity of the diagram, observe that
we have a diffeomorphism of
four-manifolds $W(G'(v))\cong W(G_{+1}(v)\cup f)$, arising by sliding the
circle corresponding to $v$ (in $G'(v)$) over the circle $e$. By the handleslide
invariance of the maps induced by cobordisms (cf.\ \cite{HolDiskFour}) and the
blow-up formula, 
we have the identification
$$\Fp{W(G'(v)),K'}(\xi)=\Fp{W(G_{+1}(v)\cup f),(K,1)}(\xi')
=\Fp{W(G_{+1}(v)),K}(\xi''),$$
where here $\xi'$ and $\xi''$ are obtained from $\xi$ by equivalences of the
Heegaard diagrams belonging to $-Y(G'(v))$, $-Y(G_{+1}(v)\cup f)$, and 
$-Y(G_{+1}(v))$, and $K'$ is the characteristic vector of the $\SpinC$ structure
over $G'(v)$ whose characteristic vector can be written as $(K,1)$
under the change-of-basis corresponding to $W(G'(v))\cong W(G_{+1}(v)\cup f)$.
Commutativity of the square now follows.
\end{proof}

We begin with a few remarks on the case of a graph with no bad points.
Recall that for such graphs, $\HFpRed(Y(G))=0$, as established in
Theorem~\ref{HolDiskSymp:thm:FloerHomology} of~\cite{HolDiskSymp}.
This follows easily from the following lemma, whose proof we include
here for the reader's convenience:

\begin{lemma}
\label{lemma:NoBadPoints}
If $G$ is a negative-definite graph with no bad vertices, then 
$H_1(Y(G);\Z)=\Rk\HFa(Y(G))$.
\end{lemma}

\begin{proof}
To show that $|H_1(Y(G);\Z)|=\Rk\HFa(Y(G))$,
one shows that both of these numbers
are additive in the sense that:
\begin{equation}
\label{eq:Additive}
N(G)=N(G_{+1}(v))+N(G-v)
\end{equation}
(provided that $G_{+1}(v)$ also has no bad points)
and they both satisfy the  normalization that 
\begin{equation}
\label{eq:Normalization}
N({\text{empty graph}})=1.
\end{equation}
(Additivity 
of $\Rk\HFa(Y(G))$ follows easily from the long exact sequence of $\HFa$;
additivity of $|H_1(Y(G);\Z)|$  is elementary.)

The equality of the two quantities now follows from induction, 
and the following observation: if $G$ is a graph with no bad vertices 
and $v$ is a leaf (i.e.\ a vertex with degree $d(v)=1$)
with multiplicity $-1$, then $G$ can be ``blown down''
to produce a graph with no bad vertices and one fewer vertex.
\end{proof}

\begin{lemma}
\label{lemma:NoBadPoints2}
Let $G$ be a graph which satisfies the inequality at each vertex $v$:
\begin{equation}
\label{ineq:StrongNoBad}
m(v)<-d(v).
\end{equation}
Then, the rank of $\Ker U\subset \Combp(G)$ is the number of 
$\SpinC$ structures over $Y$.
\end{lemma}

\begin{proof}
In view of Lemma~\ref{lemma:Dualize}, the rank of $\Ker U$ is
determined by the number of inequivalent characteristic
vectors in ${\mathcal B}_0$,
which are also not equivalent to elements in
$\Z^{>0}\times\Char(G)$. Indeed, by subtracting off
Poincar\'e duals of vertices as required, 
this is equal to the number of distinct characteristic
vectors which are not equivalent to vectors of the form
$U^n\otimes K'$ with $n>0$ and which satisfy the inequality:
\begin{equation}
\label{eq:BetterBound}
m(v)+2\leq K(v) \leq -m(v)
\end{equation}
at each vertex $v$. We call such a  characteristic vector a {\em short vector for $G$}, and let $\Short(G)$ denote the set of short vectors.  The proposition, then, is equivalent to the
statement that the number of vectors in $\Short(G)$ agrees with
the order of $H_1(Y(G);\Z)$.

To this end, let $v$ be a leaf of $G$, and $w$ be a neighbor of $v$. 
We claim that
\begin{equation}
\label{eq:ShortVectorRelation}
|\Short(G)|=-m(v) |\Short(G-v)| - |\Short(G-v-w)|.
\end{equation} 
This equation implies the lemma by induction on the number of vertices
in the graph, since $(-1)^{|G|}\det(G)$ (which counts the order of
$H_1(Y(G);\Z)$ and hence the number of $\SpinC$ structures over
$Y(G)$) clearly satisfies the same relation.

First, we claim that for each short vector $K$ for $G-v-w$, there is
some constant $m(w)+2\leq c(K)$ with the property that $(K,p)$ is a
short vector for $G-v$ for all $m(w)+2\leq p \leq c(K)$; and indeed, all
the short vectors for $G-v$ arise in this manner. Then,
Equation~\eqref{eq:ShortVectorRelation} follows from the following
claim: the set of short vectors for $G$ whose restriction to $G-v-w$
is $K$ is given by: $$\left(\bigcup_{
\begin{tiny}
{\begin{array}{c}
m(w)+2\leq p\leq c(K), \\
m(v)+2\leq i \leq -m(v)
\end{array}}
\end{tiny}} (K,p,i)\right) - (K,c(K),-m(v)).$$
To see the existence of $c(K)$ as above, we proceed as follows.  A
characteristic vector $L$ (for any negative-definite, weighted graph)
which satisfies inequalities $m(v)+2\leq L(v)\leq -m(v)$ at each
vertex $v$ is equivalent to a vector of the form $U^n\otimes L'$ with
$n>0$ if and only if there is a subset $\{v_1,\ldots,v_k\}$ of vertices
and an element $A=a_1\PD[v_1]+\cdots+a_k\PD[v_k]$ with all $a_i>0$, so
that 
\begin{equation}
\label{ineq:DiffPos}
L\cm A + A\cm A >0
\end{equation}
(see Proposition~\ref{prop:KerU} below). Moreover, it follows easily
from this same discussion that if $G$ satisfies
Inequality~\eqref{ineq:StrongNoBad} at each vertex, then we can arrange
for all the
$a_i\in \{0,1\}$.  We use this to conclude that if $K$ is a short
vector for $G-v-w$, then $L=(K,m(v)+2)$ is a short vector for $G-v$:
if not, we would have an expression $A$ as above, and, since $K$ is
short for $G-v-w$, $w\in
\{v_1,\ldots,v_k\}$, so we write $A=A_0+\PD[w]$. Let $j$ be the
number of spheres appearing in the expression for $A_0$ (with non-zero
multiplicity) which meet $w$. Then,
\begin{eqnarray*}
L\cm A + A \cm A  
&=& K\cm A_0 + A_0\cm A_0 + 2m(w)+2 + 2j.
\end{eqnarray*}
Since $j\leq d(w)<-m(w)$, and $K\in\Short(G-v-w)$, it follows
that $L\cm A + A\cm A<0$, contradicting Inequality~\eqref{ineq:DiffPos}.

The count of short vectors in $G$ with fixed restriction to $G-v-w$
(which establishes Equation~\eqref{eq:ShortVectorRelation})
proceeds similarly.
\end{proof}

Suppose $G$ is a graph with a distinguished vertex $v$.
Clearly, $Y(G)$ is obtained from $Y(G-v)$ by a single two-handle
addition, and $Y(G'(v))$ is also obtained from $Y(G)$ by a single
two-handle addition. These two-handle additions can be viewed as
cobordisms from $-Y(G)$ to $-Y(G-v)$ and $-Y(G'(v))$ to $-Y(G)$
respectively. The induced maps fit into a short exact sequence, as follows:

\begin{prop}
\label{prop:SurgerySeq}
Suppose that $G_{+1}(v)$ is a negative-definite 
plumbing diagram, and suppose that $G-v$ contains no bad points. Then,
there is a short exact sequence: 
$$
0 \to \HFp(-Y(G'(v))) \buildrel\MapOnep\over\longrightarrow
\HFp(-Y(G)) \buildrel\MapTwop\over\longrightarrow\HFp(-Y(G-v)) \to0,
$$
where the maps $\MapOnep$ and $\MapTwop$ above are induced by the two-handle additions.
\end{prop}

\begin{proof}
Theorem~\ref{HolDiskTwo:thm:GeneralSurgery} of~\cite{HolDiskTwo}
gives a long exact sequence 
\begin{gather*}
\cdots \longrightarrow\HFp(-Y(G'(v))) 
\stackrel\MapOnep{\longrightarrow}
\HFp(-Y(G))\hspace{1.5in}\\
\hspace{1in}\stackrel{\MapTwop}{\longrightarrow}
\HFp(-Y(G-v)) 
\stackrel{\MapThreep}{\longrightarrow}
\HFp(-Y(G'(v))) \longrightarrow \cdots
\end{gather*}
where the maps $\MapOnep$, $\MapTwop$, and $\MapThreep$ 
are induced by two-handle additions
(in the sense of~\cite{HolDiskFour}).

The hypothesis that $G_{+1}(v)$ is negative-definite ensures that the
two cobordisms inducing $\MapOnep$ and $\MapTwop$ above are both
negative-definite; it follows that the third cobordism is not, and
hence it induces the trivial map on $\HFinf$
(cf.\ Lemma~\ref{HolDiskFour:lemma:BTwoPlusLemma}
of~\cite{HolDiskFour}). Since $G-v$ has no bad points, it follows
easily from Lemma~\ref{lemma:NoBadPoints} that
$\HFpRed(-Y(G-v))=0$. Thus, the map $\MapThreep$ is trivial.
\end{proof}

Note that the above short exact sequence 
(together with the obvious induction on the graph)
suffices to prove Corollary~\ref{cor:EvenDegrees}.

To prove Theorem~\ref{thm:SomePlumbings}, we compare the short exact
sequence of Proposition~\ref{prop:SurgerySeq} with a corresponding 
sequence for
$\Combp$, defined as follows.

We write characteristic vectors for $G'(v)$ as triples $(K,i,\ell)$, where $K$ is the restriction
to $G-v$, and $i$ and $\ell$ denote the evaluations on
the vertices $v$ and $f$ respectively
(in particular, here $i\equiv m(v)\pmod{2}$, and $\ell\equiv 1\pmod{2}$).
Similarly, we write characteristic vectors for $G$ as pairs $(K,i)$.

Given a finitely supported map $\phi\in\Hom(\Char(G'(v)),\InjMod{0})$,
we let $$\MapOneCombp(\phi)\in\Hom(\Char(G),\InjMod{0})$$ be the map
defined by $$\langle \MapOneCombp(\phi),(K,p)\rangle
=\sum_{j=-\infty}^{+\infty} \phi(K,p,2j+1).$$ Similarly, given a
finitely supported map $\phi\in\Hom(\Char(G),\InjMod{0})$, we let
$$\MapTwoCombp(\phi)\in\Hom(\Char(G-v),\InjMod{0})$$ be the map
defined by $$\langle \MapTwoCombp(\phi),K\rangle
=\sum_{i=-\infty}^{+\infty}
\phi(K,2i+m(v)).$$

\begin{lemma}
The above formulas induce maps
\begin{eqnarray*}
\MapOneCombp\colon \Combp(G'(v))\longrightarrow \Combp(G)
&{\text{and}}&
\MapTwoCombp\colon \Combp(G)\longrightarrow \Combp(G-v).
\end{eqnarray*}
\end{lemma}

\begin{proof}
The proof is straightforward.
\end{proof}

Next, we verify that these maps fit together with the maps of
Proposition~\ref{prop:SurgerySeq} as follows:
\begin{equation}
\label{eq:CommDiagram}
\begin{CD}
\HFp(-Y(G'(v))) @>\MapOnep>>\HFp(-Y(G)) @>{\MapTwop}>> \HFp(-Y(G-v)),\\
@V{\NatTransp_{G'(v)}}VV @V{\NatTransp_{G}}VV
@V{\NatTransp_{G-v}}VV \\
\Combp(G'(v))) @>\MapOneCombp>>\Combp(G) 
@>\MapTwoCombp>> \Combp (G-v).
\end{CD}
\end{equation}

\begin{lemma}
\label{lemma:NearlyExact}
Let $G$ be a graph with the property that $G'(v)$ is
negative-definite. Then, the squares in Diagram~\eqref{eq:CommDiagram}
commute. Moreover, $\MapOneCombp$ is injective, and
$$\MapTwoCombp\circ\MapOneCombp=0.$$
\end{lemma}

\begin{proof}
Commutativity of the squares follows immediately from the naturality
of the maps under composition of cobordisms.

To verify the injectivity of $\MapOneCombp$, we proceed as follows.
Note that the grading on $\InjMod{0}$ induces a filtration in the
natural way. Specifically, if $\xi\in\InjMod{0}$ is a non-zero
element, we can 
write 
$$\xi=\xi_0+\ldots+\xi_n,$$ where $\xi_i$ denotes the
homogeneous component of $\xi$ in dimension $i$,
so that $\xi_n\neq 0$. Then, we define $\Filt(\xi)=n$. By convention, 
$\Filt(0)=-\infty$.

Now, fix some non-zero $\phi\in\Combp(G'(v))$, and let ${\mathcal K}_1\subset
\Char(G'(v))$ denote the set of characteristic vectors $K$ which
maximize $\Filt(\phi(K))$ (amongst all characteristic vectors). Such a
vector can be found since $\phi$ has finite support. Next, let
${\mathcal K}_2\subset {\mathcal K}_1$ denote the subset for which
$\langle K, v\rangle$ is maximal.  We claim that if 
$K\in {\mathcal
K}_2$, then $$\MapOneCombp(\phi)(K|G)\neq 0.$$ In fact, letting
${\mathcal E}(K|G)$ denote the set of characteristic vectors for
$G'(v)$ whose restriction to $G$ agrees with the restriction of $K$,
we claim that there is a unique vector in ${\mathcal E}(K|G)$ which
maximizes ${\Filt}\circ\phi$ (so it is $K$), and that vector satisfies
$\langle K, e\rangle = -1$. To see this, observe that if $L\in
{\mathcal E}(K|G)$ satisfies $|\langle L, e\rangle | > 1$, then
$\Filt(\phi(L))<\Filt(\phi(K))$.  This follows immediately by the
adjunction relation: if $L$ satisfies this hypothesis, then by adding
or subtracting $2\PD(e)$ (depending on the sign of $\langle L,
e\rangle$), we could find a new vector $L'$ with $U^m\cm
\phi(L')=\phi(L)$ for $m>0$,
so $\Filt(\phi(L))<\Filt(\phi(L'))\leq \Filt(\phi(K))$.
Next, we claim that $\langle K, e\rangle=-1$, for if $\langle K,
e\rangle=+1$, then by adding $2\PD(e)$, we would obtain a new
characteristic vector $K'$ with $\phi(K)=\phi(K')$, but with $$\langle
K', v\rangle = \langle K,v\rangle+2,$$ violating the hypothesis that
$K\in {\mathcal K}_2$. This completes the proof that $K$ is the
unique vector in ${\mathcal E}(K|G)$
which maximizes  $\Filt\circ\phi$, so it follows immediately that
$\MapOneCombp(\phi)(K|G)\neq 0$.

To verify that $\MapTwoCombp\circ \MapOneCombp=0$, we proceed as
follows.  It is an easy consequence of the adjunction relation that if
$K\in\Char(G-v)$ is any characteristic vector, then for any $i\geq 0$, we
have that
\begin{eqnarray*}
\phi(K,p,2i+1)&=& U^{\frac{i(i+1)}{2}}\cm \phi(K,p+2i+2,-1), \\
\phi(K,p,-(2i+1))&=& U^{\frac{i(i+1)}{2}}\cm \phi(K,p-2i,-1),
\end{eqnarray*}
where $(K,i,j)\in\Char(G'(v))$ is the characteristic vector whose restriction
to $G-v$ is $K$, and its values on $v$ and $e$ are $i$ and $j$ respectively.

Thus in the double sum
\begin{eqnarray*}
\langle (\MapTwoCombp\circ\MapOneCombp)(\phi), K\rangle 
&=& \sum_{i=-\infty}^{+\infty}\sum_{j=-\infty}^{+\infty}
\phi(K,2i+m(v),2j+1),
\end{eqnarray*}
the terms cancel in pairs.
\end{proof}

\begin{prop}
\label{prop:NoBadPointsVStrong}
Suppose $G$ is a plumbing diagram satisfying
Inequality~\eqref{ineq:StrongNoBad} at each vertex, then $\NatTransp$
induces an isomorphism $$\NatTransp\colon
\HFp(-Y(G))\stackrel{\cong}{\longrightarrow} \Combp(G).$$
\end{prop}

\begin{proof}
With a slight abuse of notation, we let 
$$\Comba(G(v))=\Ker U\subset \Combp(G(v))$$
throughout this proof.

We prove the result by induction on the graph.
The basic cases where the graph is empty (so there is only one characteristic
vector, the zero vector) is obvious, as is the case where the graph 
has a single vertex labelled with multiplicity $-1$. 

We now prove the result by induction on the number of vertices in the
graph, and then a sub-induction on $-m(v)$, where $v$ is a leaf in the
graph.  For the sub-induction, we allow $m(v)=-1$. This case is
handled by the inductive hypothesis on the number of vertices, since
$-Y(G)$ is equivalent to a plumbing on a graph with fewer vertices and
no bad points,
and observing that the identification from Proposition~\ref{prop:BlowDown}
is natural under $T^+$.

For the sub-inductive step on $-m(v)$, consider the following analogue
of Diagram~\eqref{eq:CommDiagram} (which also commutes, according to
Lemma~\ref{lemma:NearlyExact}): $$
\begin{CD}
0 \longrightarrow\HFa(-Y(G'(v))) @>\MapOnea>>\HFa(-Y(G)) @>{\MapTwoa}>> \HFa(-Y(G-v)) \longrightarrow0,\\
 @V{\NatTransa_{G'(v)}}VV @V{\NatTransa_{G}}VV
@V{{\NatTransa}_{G-v}}VV \\
 \Comba(G'(v))) @>\MapOneComba>>\Comba(G) 
@>\MapTwoComba>> \Comba (G-v).
\end{CD}
$$ 
Here, the top row is exact according to Proposition~\ref{prop:SurgerySeq}.
It follows by the inductive hypotheses, that both
$\NatTransa_{G'(v)}$ and $\NatTransa_{G-v}$ are isomorphisms. For
$\NatTransa_{G-v}$, this is obvious, while for $\NatTransa_{G'(v)}$ we
use Proposition~\ref{prop:BlowDown}. It follows immediately that
$\MapTwoComba$ is surjective.  Moreover, since
$\MapTwoComba\circ\MapOneComba=0$, and the rank of $\Comba(G)$ is the
sum of the ranks of $\Comba(G'(v))$ and $\Comba(G-v)$ (according to
Lemma~\ref{lemma:NoBadPoints2}), it follows from an easy count of
ranks that the kernel of $\MapTwoComba$ is the image of
$\MapOneComba$. Now, by the five-lemma, it follows that $\NatTransa_G$
is an isomorphism.

We now consider Diagram~\eqref{eq:CommDiagram}, where $v$ is a leaf.
Observe that $\MapOnep$ is injective while $\MapTwop$ is surjective
according to Proposition~\ref{prop:SurgerySeq}.
Again, by induction we have that $\NatTransp_{G'(v)}$ and
$\NatTransp_{G-v}$ are isomorphisms, so that $\MapOneCombp$ is injective
and $\MapTwoCombp$ is surjective. Now, exactness in the middle follows
easily from the fact that it holds on the level of $\HFa$, together
with the fact that $\MapTwoCombp\circ \MapOneCombp=0$. As before, we can
now use the five-lemma to establish the desired isomorphism.
\end{proof}

\begin{prop}
\label{prop:NoBadPointsStrong}
Suppose $G$ is a negative-definite plumbing diagram with no bad
points, then $\NatTransp$ induces an isomorphism $$\NatTransp\colon
\HFp(-Y(G))\stackrel{\cong}{\longrightarrow} \Combp(G).$$
\end{prop}

\begin{proof}
This is proved by induction on the number of vertices $G$ with
$d(v)=-m(v)$. The case where there is no vertex with $d(v)=-m(v)$, is
Proposition~\ref{prop:NoBadPointsVStrong}. For the inductive step, we
consider Diagram~\eqref{eq:CommDiagram} again (together with
Proposition~\ref{prop:SurgerySeq}), where $v\in G$ is a vertex (not
necessarily a leaf) with $d(v)=-m(v)-1$. This time induction tells us
that $\NatTransp_{G-v}$ and $\NatTransp_G$ are isomorphisms. This,
together with Lemma~\ref{lemma:NearlyExact}, is sufficient (after a
straightforward diagram chase) to allow us to conclude that
$\NatTransp_{G'(v)}$ is an isomorphism.
\end{proof}

\vskip.2cm
\proof[Proof of Theorem~\ref{thm:SomePlumbings}]
This follows in from Proposition~\ref{prop:NoBadPointsStrong} in a
manner analogous to how that proposition follows from
Proposition~\ref{prop:NoBadPointsVStrong}.  Again, we consider
Diagram~\eqref{eq:CommDiagram}, now choosing $v$ to be the bad
vertex. This time, induction and
Proposition~\ref{prop:NoBadPointsStrong} tells us that
$\NatTransp_{G-v}$ and $\NatTransp_G$ are isomorphisms. Again, we
conclude (with the help of Lemma~\ref{lemma:NearlyExact}) that
$\NatTransp_{G'(v)}$ is an isomorphism.
\endproof

\proof[Proof of Theorem~\ref{thm:PushedFurther}]
We proceed as in the above proof.  Let $v$ be one of the bad vertices
in $G$. We would like to prove the result by descending induction on
$-m(v)$. In this case, however, Proposition~\ref{prop:SurgerySeq} is
no longer available to us since $G-v$ has a bad vertex; but it is the
case that $\HFp_{\odd}(-Y(G-v))=0$, in view of
Theorem~\ref{thm:SomePlumbings} (or, more precisely,
Corollary~\ref{cor:EvenDegrees}).  Thus, we have the diagram: $$
\begin{CD}
0 @>>>\HFp_{\ev}(-Y(G'(v))) @>\MapOnep>>\HFp_{\ev}(-Y(G)) @>{\MapTwop}>> \HFp_{\ev}(-Y(G-v))\\
&& @V{\NatTransp_{G'(v)}}VV @V{\NatTransp_{G}}V{\cong}V
@V{\NatTransp_{G-v}}V{\cong}V \\
0@>>>  \Combp(G'(v))) @>\MapOneCombp>>\Combp(G) 
@>\MapTwoCombp>> \Combp (G-v),
\end{CD}
$$
where the maps are indicated as isomorphisms when it follows from induction. This
diagram forces $\NatTransp_{G'(v)}$ to be an isomorphism, as well.
\endproof

\proof[Proof of Corollary~\ref{cor:CalcD}]
Fix a $\SpinC$ structure $\spinct$ over $Y$, and let
$K_0$ be a characteristic vector in $\Char_{\spinct}(G)$ 
for which $K^2$ is maximal. We define a sequence 
of elements $\Phi_N\in\Combp(G,\spinct)$ by
$$\Phi_N(K)=U^{\left(\frac{K_0^2-K^2}{8}\right)-N}\in\InjMod{0}.$$
As usual, if the exponent of $U$ here is positive, the corresponding
element of $\InjMod{0}$ is zero.
so, in particular, 
$$\Phi_0(K)=\left\{\begin{array}{ll}
1 & {\text{if $K^2$ is maximal in $\Char_{\spinct}(G)$}}\\
0, & {\text{otherwise}}
\end{array}\right.$$
so 
$$\deg(\Phi_0)=-\left(\frac{K_0^2+|G|}{4}\right).$$
Clearly, 
$U\cm \Phi_{N+1}=\Phi_{N}$,  and $U\cm \Phi_0=0$.
Thus, by Theorem~\ref{thm:SomePlumbings} (and 
Theorem~\ref{thm:PushedFurther}, in the case where there are two
bad points), it follows that 
$\deg(\Phi_0)=d(-Y(G),\spinct)$. Since $d(-Y(G),\spinct)=-d(Y(G),\spinct)$,
the corollary follows.
\endproof

\subsection{Theorem~\ref{intro:SomePlumbings} over $\Z$}
\label{subsec:Signs}

Strictly speaking, when working over $\Z$, the map associated to a
cobordism as defined in~\cite{HolDiskFour} does not have a canonical
sign. Thus, it might appear that $T^+$ is defined only as a
map $$T^+\colon \HFp(-Y(G),\spinct)\longrightarrow
\Hom(\Char_\spinct(G),\InjMod{0}/\pm 1).$$ 
In fact, we can actually specify
a map $$T^+\colon \HFp(-Y(G),\spinct)\longrightarrow
\Hom(\Char_\spinct(G),\InjMod{0}),$$ which is well-defined 
up to an overall $\pm 1$ sign, which we pin down with some additional
data.  Specifically, fix a single $\SpinC$ structure
$\spinc\in\SpinC(W(G))$ whose restriction to $-Y(G)$ is
$\spinct$. Since $W(G)$ is a negative-definite cobordism between
rational
homology spheres, the induced map $$\Finf{W,\spinc}\colon
\HFinf(-Y(G),\spinct)\longrightarrow
\HFinf(S^3)$$ is an isomorphism (cf.\ Proposition~\ref{AbsGraded:prop:NegSurgery} of~\cite{AbsGraded}), and hence determined up to an overall
$\pm 1$. Now, for each other $\SpinC$ structure $\spinc'$, we choose
orientation conventions so that the induced isomorphism
$\Finf{W,\spinc'}$ agrees with $U^m\cm \Finf{W,\spinc}$ (where here,
of course, $m=\frac{c_1(\spinc)^2-c_1(\spinc')}{8}$).  This fixes
signs for all the maps $\Fp{W(G),\spinc'}$.

With this sign in place, we see that $T^+$ induces a map to
$\Combp(G,\spinct)$ with $\Z$ coefficients (i.e.\ Proposition~\ref{prop:SetUp} works over $\Z$),
which is uniquely determined up to one overall sign.

Before defining $\MapOneCombp$ and $\MapTwoCombp$, we pause for a discussion
of their geometric counterparts $\MapOnep$ and $\MapTwop$ appearing in the 
long exact sequence for $\HFp$.

Recall that there was considerable leeway in the orientation
conventions used in defining the maps in the surgery long exact sequence for
$\HFp$, see~\cite{HolDiskTwo}. Indeed, the maps $A^+$ and $B^+$ are
defined as sums of maps induced by cobordisms, and any orientation
convention was allowed provided that the composite map (on the chain
level) is chain homotopic to $0$. More concretely, this can be
stated (using notation from Lemma~\ref{lemma:NearlyExact}) as follows.
Let $W_1$ be the cobordism from $-Y(G'(v))$ to $-Y(G)$ and $W_2$ be
the cobordism from $-Y(G)$ to $-Y(G-v)$, and observe that the
composite cobordism can be identified as a blow-up $W_1\cup
W_2=W_3\#\mCP$ (with exceptional sphere $e$), then we can let $A^+$ and $B^+$ be the maps on
homology induced by chain maps:
\begin{eqnarray*}
a^+=\sum_{\spinc_1\in\SpinC(W_1)}\alpha(\spinc_1)\cm
\fp{W_1,\spinc_1} &{\text{and}}&
b^+=\sum_{\spinc_2\in\SpinC(W_2)}\beta(\spinc_2)\cm
\fp{W_2,\spinc_2},
\end{eqnarray*}
where $\fp{}$ denotes the chain map induced by the cobordism
(with the canonical orientation convention, since both 
are negative-definite cobordisms in our case),
and 
\begin{eqnarray*}
\alpha\colon \SpinC(W_1)\longrightarrow \{\pm 1\},
&&
\beta\colon \SpinC(W_2)\longrightarrow \{\pm 1\}
\end{eqnarray*}
are maps satisfying the constraint that if $\spinc, \spinc'\in
\SpinC(W_1\# W_2)$ are any two $\SpinC$ structures which agree over
$W_3$, and $$\langle c_1(\spinc), e\rangle = -\langle c_1(\spinc'), e
\rangle,$$ then $$\alpha(\spinc|_{W_1})\cm \beta(\spinc|_{W_2})=
-\alpha(\spinc'|_{W_1})\cm \beta(\spinc'|_{W_2}).$$ For example, we
can choose $\beta\equiv 1$, and $\alpha$ as follows. Let $\PD[e]\in
W_1\# W_2$, and let $\epsilon=\PD[E]|_{W_1}$. For each
$\Z\epsilon$-orbit in $\SpinC(W_1)$, fix an initial $\SpinC$
structure $\spinc_0$ over $W_1$, and let $\alpha(\spinc_0)=1$. Then, if
$\spinc-\spinc_0= m \epsilon$, we let $\alpha(K)=(-1)^m$.

We now define the maps fitting into the short exact sequence:
\begin{eqnarray*}
\MapOneCombp\colon \Combp(G'(v))\longrightarrow \Combp(G),
&{\text{and}}&
\MapTwoCombp\colon \Combp(G)\longrightarrow \Combp(G-v).
\end{eqnarray*}
These maps are defined by the  formulas:
\begin{eqnarray*}
\langle \MapOneCombp(\phi),(K,p)\rangle
&=&\sum_{j=-\infty}^{+\infty} \alpha(K,p,2j+1)\phi(K,p,2j+1), \\
\langle \MapTwoCombp(\phi),K\rangle
&=&\sum_{i=-\infty}^{+\infty}
\beta(K,2i+m(v))\phi(K,2i+m(v)),
\end{eqnarray*}
where here $\alpha(K,p,2j+1)$ denotes $\alpha$ applied to the restriction
to $W_1$ of the $\SpinC$ structure over $W(G'(v))$ whose first Chern class
is $(K,p,2j+1)$, with the similar shorthand for $\beta$.

With these sign conventions in place, we claim that 
the analogue of Lemma~\ref{lemma:NearlyExact} now holds.
(Where the statement about commutative squares is weakened to squares
which commute, up to sign.)
Indeed, all the proofs from the last subsection readily
adapt now to prove both Theorems~\ref{thm:SomePlumbings} and
\ref{thm:PushedFurther} over $\Z$.

\section{Calculations}
\label{sec:Examples}

\subsection{An algorithm for determining the rank of $\Ker U$}

The group $\Combp(G)$ can be determined from the combinatorics of the
plumbing diagram. In fact, Lemma~\ref{lemma:Dualize} gives us a finite
model for $\Combp(G)$ (at least, the subset of $\Ker U^{n+1}$, for
arbitrary $n$).  We give here a more practical algorithm for
calculating $\Ker U$ or, more precisely, its dual space.

Fix a characteristic vector $K$ satisfying inequality
\begin{equation}
\label{eq:PartBox}
m(v)+2 \leq \langle K, v\rangle \leq -m(v).
\end{equation}
Now, successively apply the following algorithm to find a path of vectors
$(K=K_0,K_1,\ldots,K_n)$ so that the $K_i$ for $i<n$ satisfy the bounds
\begin{eqnarray}
\label{eq:InBigBox}
|\langle K_i,v\rangle|\leq -m(v)
\end{eqnarray} for all vertices $v$.
Given $K_{i}$, choose any vertex $v_{i+1}$ with
\begin{eqnarray}
\label{eq:Induct}
\langle K_i,v_{i+1}\rangle = -m(v_{i+1}),
&{\text{then let}}&
K_{i+1}=K_i+2\PD[v_{i+1}].
\end{eqnarray}
This algorithm can terminate in one of two ways:
either
\begin{itemize}
\item the final vector $L=K_n$ satisfies the inequality,
\begin{equation}
\label{eq:OtherPartBox}
m(v) \leq \langle L, v\rangle \leq -m(v)-2
\end{equation}
at each vertex $v$ or
\item there is some vertex $v$ for which 
\begin{equation}
\label{eq:OtherTermination}
\langle K_{n},v \rangle > -m(v).
\end{equation}
\end{itemize}

To calculate the rank of $\Ker U$ in the examples in this paper, we
use the following claim: the equivalence classes in $\DCombp(G)$ which
have no representative of the form $U^m\otimes K'$ with $m>0$ are in
one-to-one correspondence with initial vectors $K$ satisfying
Inequality~\eqref{eq:PartBox} for which the algorithm above terminates
in a characteristic vector $L$ satisfying
Inequality~\eqref{eq:OtherPartBox}.  The purpose of
Proposition~\ref{prop:KerU} is to establish this claim. 

\begin{defn}
\label{def:FullPath}
We call a sequence of characteristic vectors 
$(K=K_0,\ldots,K_n$ $=L)$ obtained from the above algorithm
a {\em full path}; i.e.\
$K=K_0$ satisfies Inequality~\eqref{eq:PartBox}, $K_{i+1}$ is obtained from $K_i$
by Equation~\eqref{eq:Induct}, and the final vector $L=K_n$ satisfies either
Inequality~\eqref{eq:OtherPartBox} or~\eqref{eq:OtherTermination}.
\end{defn}

\begin{prop}
\label{prop:KerU}
Fix an equivalence class in $\DCombp(G)$ which contains no
representatives of the form $U^m\otimes K'$ for $m>0$. Each such equivalence
class
has a unique representative $K$ satisfying Equation~\eqref{eq:PartBox}.
Indeed, a characteristic vector $K$ satisfying these bounds
is inequivalent to an element of the form $U^m\otimes K'$ (with $m>0$)
if and only if we can find a full path $$(K=K_0,K_1,\ldots,K_n=L)$$ terminating
with a characteristic vector $L=K_n$ which satisfies Inequality~\eqref{eq:OtherPartBox}
for each vertex $v$.
\end{prop}

\begin{proof}
Let $M$ be a characteristic vector which is not equivalent to
$U^m\otimes K'$ for $m>0$.  We find a full path using the above
algorithm. Specifically we let $L_0=M$, and then for each $j\geq 0$,
and extend $L_0$ to a sequence $L_0,\ldots,L_{n_+}$ by letting $v_{j+1}$
be a vertex for which $$\langle L_j, v_{j+1}\rangle = -m(v_{j+1}),$$
and then letting $L_{j+1}=L_j+2\PD[v_{j+1}]$.  Clearly, in this
sequence, each element satisfies $$m(v)\leq \langle L_j, v\rangle
\leq  -m(v)$$ (for otherwise, $L_j$ would be equivalent to an element of
$\Z^{>0}\times\Char(G)$). The sequence is finite (since the elements
of the sequence are all distinct), so it must terminate with $L_{n_+}$
satisfying Inequality~\eqref{eq:OtherPartBox}. 
In the same way, we can extend back from $L_0$ to obtain a sequence
$(L_0,L_{-1},\ldots,L_{n_-})$ by the rule that
if there is a vertex $v$ for which $\langle L_j,v \rangle = m(v)$,
then $$L_{j-1}=L_j-2\PD[v].$$ This sequence must terminate with $L_{n_-}$ satisfying
Inequality~\eqref{eq:PartBox}. Thus, 
$$(L_{n_-},L_{n_-+1},\ldots,L_0=M,L_1,\ldots,L_{n_+})$$
is a full path in the 
sense of Definition~\ref{def:FullPath}.
In particular, $L_{n_-}$ is the representative $K$. 

We argue that if $M$ is a vector in a full path
$(K=K_0,\ldots,K_n=L)$, so that $L$ satisfies
Inequality~\eqref{eq:OtherPartBox}, then $L$ is uniquely determined by
$M$ (i.e.\ independent of the particular sequence). In fact, if
$\{v_1,\ldots,v_\ell\}$ are vertices with $\langle M,v_i\rangle =
-m(v_i)$, then $L$ must be obtained from $M$ by adding
$2\PD[v_1]+\cdots+2\PD[v_\ell]$ (so that we can achieve $\langle
L,v_i\rangle<-m(v_i)$), and then adding some additional
vertices. Thus, 
$$M+2\PD[v_1]+\cdots+2\PD[v_\ell]$$
lies on a full path with
the same endpoint $L$. By induction on the minimal distance of $M$ to
its endpoint on a full path, we have the uniqueness of the final point
$L$. 

Next, we argue that if $M$ and $M'$ are two characteristic vectors
which are equivalent to one another, and $M\not\sim U^m\otimes K'$ for
$m>0$, then the endpoint of any full path through $M$ agrees with the
endpoint of a full path through $M'$. This is clear if $M'=M\pm
2\PD[v]$: we can find a full path which passes through both $M$ and
$M'$. More generally, if $M\sim M'$, we can get from $M$ to $M'$ by a
finite number of additions or subtractions of $2\PD[v]$ for vertices
so as to leave the square unchanged (i.e.\ $M_{i+1}$ is obtained from
$M_i$ by $M_{i+1}=M_i- 2\epsilon_{i+1} \PD[v_{i+1}]$ where $v_{i+1}$
is a vertex which satisfies $\langle M_i,v_{i+1}\rangle =
\epsilon_{i+1} m(v_{i+1})$ for $\epsilon_{i+1}=\pm 1$). The assertion
then follows by an easy induction on the number of such operations.

Turning this around, we also see that the initial point of a full path 
is uniquely determined by the equivalence class of the characteristic 
vectors lying in it. This gives the uniqueness statement claimed in the 
proposition.

Finally, we argue that if $K$ is a vector satisfying
Inequality~\eqref{eq:PartBox}, but $K\sim U^m \otimes K'$ for $m>0$,
then there is no full path connecting $K$ to another characteristic
vector $L$ satisfying Inequality~\eqref{eq:OtherPartBox}. To see this
suppose that $K\sim U^m\otimes K'$, we can find some sequence $$M_0=K,
M_1,\ldots,M_\ell=M$$ and signs $\epsilon_i\in \{\pm 1\}$ with 
$$M_{i+1}=M_i-2\epsilon_{i+1}
\PD[u_{i+1}]$$ and 
$$\langle M_i,[u_{i+1}]\rangle = \epsilon_{i+1} m(u_{i+1}).$$ where each
$M_i$ satisfies (for each vertex $v$) $$|\langle M_i,[v]\rangle |\leq
-m(v),$$ but there is some vertex $w$ with $|\langle M, [w]\rangle|>
-m(v)$. We claim that
this sequence can be shortened so that all the $\epsilon_i$ are
positive.  To see this, observe first that $\epsilon_1=+1$.
Now,
consider the smallest integer $i$ with $\epsilon_{i+1}=-1$, so that 
$M_i\cm u_{i+1}=m(u_{i+1})$. It follows easily that 
$u_{i+1}\in\{u_1,\ldots,u_{i}\}$. If we
let $j$ be the last integer in $[1,\ldots,i]$ with $u_j=u_{i+1}$, then it
is easy to see that for each $k\in [j+1,\ldots,i]$, $u_k\cm u_{i+1}=0$. It
follows immediately that we delete the occurance of the $j^{th}$ and
$(i+1)^{st}$ vertices from $(u_1,\ldots,u_{n})$ to construct a new, shorter
sequence on characteristic vectors satisfying the same properties as
the $M_i$ (in particular, connecting the same two endpoints), only
with one fewer occurance of the sign $\epsilon_{j}=-1$. Proceeding in
this manner, we end up with a sequence with all
$\epsilon_j=+1$. 

Next, suppose that there is a full path connecting $K$ as above to
$L$.  We claim now that, after possibly reordering, $(u_1,\ldots,u_\ell)$
is a subsequence of the vertices $(v_1,\ldots,v_n)$ belonging to the
hypothesized full path connecting $K$ to $L$. It is easy to see then
that we can extend the original sequence $(K=M_0,\ldots,M_\ell)$ to a
sequence $K=M_0,\ldots,M_n=L$ (using a reordering $(w_1,\ldots,w_n)$ of
$(v_1,\ldots,v_n)$) so that $$M_{i+1}=M_i+2\PD[w_{i+1}]$$ and $\langle
M_{i},[w_{i+1}]\rangle \geq -m(w_{i+1})$.  This forces $K\sim L\sim 
U^m\otimes L$ for $m>0$. But it is impossible for $U^m\otimes L\sim L$.
\end{proof}

\subsection{Examples}
\label{subsec:Examples}

We illustrate the algorithm described above by  calculating $\HFp(Y)$ for 
certain Brieskorn spheres $Y$.

\medskip
{\bf{Notational Conventions}}\qua  In describing graded $\Z[U]$-modules, we adopt
the following conventions.  $\InjMod{k}$ will denote the graded
$\Z[U]$-module which is isomorphic as a relatively graded
$\Z[U]$-module to $\HFp(S^3)$, but whose bottom-most non-zero
homogeneus element has degree $k$. Also, $\Z_{(k)}$ will denote the
$\Z[U]$
module $\Z[U]/U\Z[U]$, graded so that is supported in degree $k$.

For example, with this notation, 
$$\HFp(S^2\times S^1)\cong {\mathcal T}^+_{-1/2}\oplus {\mathcal T}^+_{1/2}$$

{\bf{The Poincar\'e homology sphere}}\qua  Consider the Poincar\'e
homology sphere $Y=\Sigma(2,3,5)$. This can be realized as the
boundary of the plumbing of spheres specified by the negative-definite
$E_8$ Dynkin diagram.

We claim that the techniques of the present paper can be used to verify that
$$\HFp(-\Sigma(2,3,5))\cong {\mathcal T}^+_{-2},$$
compare Section~\ref{AbsGraded:sec:SampleCalculations} of~\cite{AbsGraded}.

We claim that the only full path connecting vectors $K$ and $L$ as in
Proposition~\ref{prop:KerU} is the path consisting of the single
characteristic vector $K=L=0$. 

Specifically, we consider the $256$ possible initial characteristic
vectors $K$ as in Proposition~\ref{prop:KerU}, i.e.\ 
$$\langle K,v_i\rangle
\in \{0,2\}.$$
It is easy to see that if $\langle K,v_i\rangle = 2$ for at least two vertices, then
the algorithm given above terminates with a characteristic vector $L$ satisfying
Inequality~\eqref{eq:OtherTermination}: i.e.\ $K\sim U\otimes K'$.

It remains then to rule out eight remaining cases where there is only
one vertex on which $K$ does not vanish. Ordering the vertices in the
$E_8$ diagram so that $v_1$ is the central node (with degree three),
and $(v_1,v_2)$, $(v_1,v_3,v_4)$, and $(v_1,v_5,v_6,v_7,v_8)$ are
three connected segments, we write characteristic vectors as tuples
$$(\langle K,v_1\rangle, \ldots \langle K,v_8\rangle).$$ We include here 
one of these eight cases -- exhibiting a full path from $K_0=(0,0,0,0,0,0,0,2)$
to a vector $K_n$ with $\langle K_n,v\rangle = 4$ for some vertex $v$ -- leaving
the remaining seven cases to the reader.
$$
\begin{array}{lll}
\big\{(2, 0, 0, 0, 0, 0, 0, 0), &
(-2, 2, 2, 0, 2, 0, 0, 0), &
(0, -2, 2, 0, 2, 0, 0, 0), \\
(2, -2, -2, 2, 2, 0, 0, 0), &
(-2, 0, 0, 2, 4, 0, 0, 0)\big\} &
\end{array}
$$

\noindent{\bf{The Brieskorn sphere $\Sigma(2,3,7)$}}\qua
We give here another calculation showing that 
$$\HFp(-\Sigma(2,3,7))\cong {\mathcal T}^+_0\oplus \Z_{(0)}$$
(compare~\cite{AbsGraded}).

This homology sphere is realized as the boundary of 
a plumbing diagram with a central node $v_1$ of
square $-1$, and three more spheres $v_2$, $v_3$, and $v_4$ of squares
$-2$, $-3$, and $-7$ respectively.  The vectors $(1, 0, -1, -3)$, $(1,0, -1, -5)$ 
are the only two vectors satisfying Inequality~\eqref{eq:PartBox} which begin a full path
ending in a characteristic vector as in Inequality~\eqref{eq:OtherPartBox}. For convenience,
we include the full path starting at $(1,0,-1,-5)$:
$$
\begin{array}{llll}
\{ (1,0,-1,-5), & (-1,2,1,-3), & (1,-2,1,-3), & (-1,0,3,-1), \\
   (1,0,-3,-1), & (-1,2,-1,1), & (1,-2,-1,1), & (-1,0,1,3) \}.
\end{array}
$$
(The other full path  is obtained by multiplying all above vectors by $-1$ and reversing the order.)
Now, we claim that $U\otimes (1,0,-1,-3)\sim U\otimes (1,0,-1,-5)$.
In fact, it is straightforward to verify that:
$$U\otimes (-1,0,1,5)\sim (1,0,1,-9)\sim (-1,0,5,-5)\sim U\otimes (1,0,-1,-5).$$
Here, we have broken the equivalence up so that when we write $K\sim K'$,
we mean that $K'$ is obtained from $K$ by applying the algorithm for constructing
a full path.

It is interesting to note that it follows from the above calculations 
that the conjugation action, which in general gives a an
involution on $\HFp(Y)$, in the present case
permutes the two zero-dimensional generators.
Observe also that the renormalized length ($\frac{K^2+|G|}{4}$)
of both vectors is $0$.

\begin{figure}[ht!]
\cl{\epsfxsize 1.7in\epsfbox{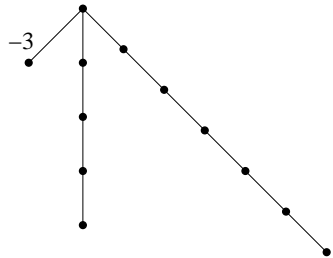}}
\caption{\label{fig:Brieskorn}
{\bf{Plumbing description of $\Sigma(3,5,7)$}}\qua
Here, the unlabeled vertices have multiplicity $-2$.}
\end{figure}

\medskip
\noindent{\bf{The Brieskorn sphere $\Sigma(3,5,7)$}}\qua
We claim that
$$
\HFp(-\Sigma(3,5,7))\cong 
{\mathcal T}^+_{-2}\oplus\Z_{(-2)}\oplus \Z_{(0)} \oplus \Z_{(0)}.
$$ We can realize $\Sigma(3,5,7)$ as the boundary of a
negative-definite plumbing of spheres, as in Figure~\ref{fig:Brieskorn}.
Unlabeled vertices
all have multiplicity $-2$. 
We order the vertices so that the central node comes first,
the $-3$-sphere second, then the four vertices on the next chain
(ordered so that the length to the central node is increasing)
and finally, the six vertices on the final chain (ordered in the same manner).

We claim that there are exactly four characteristic vectors satisfying 
Inequality~\eqref{eq:PartBox} which can be completed to
full paths terminating in characteristic vectors
satisfying Inequality~\eqref{eq:OtherPartBox}, and these are the 
vectors:
\begin{eqnarray*}
K_1&=& (0, -1, 0, 0, 0, 0, 0, 0, 0, 0, 0, 0) \\
K_2 &=& (0, 1, 0, 0, 0, 0, 0, 0, 0, 0, 0, 0) \\ 
K_3&=& (0, 1, 0, 0, 0, 0, 0, 0, 0, 0, 0, -2) \\
K_4&=& (0, 1, 0, 0, 0, -2, 0, 0, 0, 0, 0, 0)
\end{eqnarray*}
It is straightforward to verify that
\begin{eqnarray*}
U\otimes K_3 &\sim &
(2, -5, 0, 0, 0, 0, 0, 0, 0, 0, 0, -2) \\
&\sim & 
(0, -1, 0, 0, 0, 0, 0, 0, -2, 4, 0, -2) \\
&\sim & U\otimes 
(0, -1, 0, 0, 0, 0, 0, 0, 0, 0, 2, -2)
\\
&\sim & U\otimes (0, -1, 0, 0, 0, 0, 0, 0, 0, 0, 0, 2) \\
&\sim & U\otimes K_4.
\end{eqnarray*}
Here, as before, we break the equivalence up into simpler steps,
writing $K\sim K'$ if $K'$ is obtained from $K$ by applying the
algorithm for constructing a full path.

Also we have that
\begin{eqnarray*}
U\otimes K_1&\sim& (-2, 5, 0, 0, 0, 0, 0, 0, 0, 0, 0, 0) \\
&\sim& (0, 1, 0, 0, 0, 0, 0, 0, 2, -4, 0, 0) \\
&\sim& U\otimes (0,1,0, 0, 0, 0, 0, 0, 0, 0, -2, 0) \\
&\sim& U\otimes (0,-1,0, 0, 0, 0, 0, 0, 0, 2, 0, -4) \\
&\sim& U^2\otimes (0,-1,0, 0, 0, 0, 0, 0, 0, 2, -2,0) \\
&\sim& U^2\otimes K_4.
\end{eqnarray*}
A similar calculation shows that
$$U\otimes K_2 \sim U^2\otimes K_3\sim U^2\otimes K_4.$$
The result then follows.

\section{The Floer homology of $\Sigma_2\times S^1$}
\label{sec:GenusTwo}

As an application of the calculations for plumbing diagrams, we
calculate $\HFp$ of the product of a genus two surface with the
circle. Along the way, we also calculate $\HFp$ for certain other
genus two fiber-bundles over the circle.

\begin{figure}[ht!]
\relabelbox\small
\cl{\epsfxsize 2.5in\epsfbox{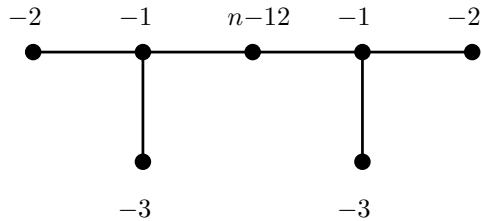}}
\relabel {-2}{${-}2$}
\relabel {-2a}{${-}2$}
\adjustrelabel <-10pt,0pt> {-1a}{${-}1$}
\relabel {-1}{${-}1$}
\relabel {-3}{${-}3$}
\relabel {-3a}{${-}3$}
\relabel {n}{$n{-}12$}
\endrelabelbox
\caption{\label{fig:Tref2}
{\bf{Plumbing description of a connected sum of two trefoils}}\qua Here, for an arbitrary
integer $n$, we have a description of $Y_n=S^3_n(T_2)$ as a plumbing diagram.}
\end{figure}

\begin{figure}[ht!]
\cl{\epsfxsize 1.35in\epsfbox{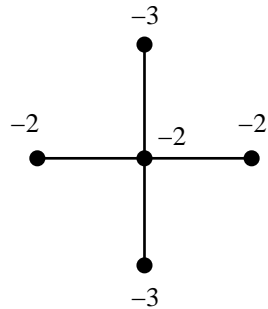}}
\caption{\label{fig:Twelve}
{\bf{Plumbing description for $Y_{12}$}}}  
\end{figure}

Let $T_2\subset S^3$ denote the connected sum of two copies of the
right-handed trefoil $T$. Now, if $Y_n=S^3_n(T_2)$ denote the
three-manifold obtained by $+n$-surgery on $S^3$ along $T_2$, then
this manifold can be realized by plumbing along the tree pictured in
Figure~\ref{fig:Tref2}.  Note that this graph has at least two bad
vertices. However, for $n=+12$, after a handleslide followed by a
handle cancellation, we obtain an alternate description of $Y_{12}$ as
the Seifert fibered space whose plumbing diagram is pictured in
Figure~\ref{fig:Twelve}.

Note that the realization of $Y_n$ as surgery on a knot gives a
correspondence $$Q\colon \Zmod{n}\longrightarrow \SpinC(Y_n).$$ We
adopt here the conventions for the integral surgery long exact
sequence, cf.\ Theorem~\ref{HolDiskTwo:thm:ExactP}
of~\cite{HolDiskTwo}.  According to these conventions
(cf.\ Lemma~\ref{AbsGraded:lemma:IdentifyQ} of~\cite{AbsGraded}), if $W$
denotes the cobordism from $S^3$ to $Y_n$ obtained by attaching the
two-handle, and $[F]\in H_2(W;\Z)$ is a generator, then $Q(0)$ is the
$\SpinC$ structure over $Y_n$ which has an extension $\spinc$ over $W$
with $$\langle c_1(\spinc),[F]\rangle = n.$$

\begin{lemma}
\label{lemma:FirstY12Calc}
Let $Q(0)$ be the $\SpinC$ structure over $Y_{12}$ as above on the
zero-surgery on the double-trefoil.  Then,
$$\HFp(-Y_{12},Q(0))\cong \Z_{(-3/4)}\oplus \InjMod{-3/4}.$$
Equivalently,
$$\HFp(Y_{12},Q(0))\cong \Z_{(-1/4)}\oplus \InjMod{3/4}.$$
\end{lemma}

\begin{proof}
Consider the plumbing diagram $G$ for $Y_{12}$ in
Figure~\ref{fig:Twelve}.  We argue that $\Ker U\subset
\Combp(G,Q(0))$ is two-dimensional.

We order the spheres $S_1,S_2,S_3,S_4,S_5$, so that $S_1$ is the
central sphere, and $S_2$ and $S_3$ are the other two two-spheres with
square $-2$. We identify characteristic vectors as quintuples,
according to the values on $S_1,\ldots,S_5$.

We use Lemma~\ref{lemma:Dualize} and Proposition~\ref{prop:KerU}
to calculate
$\Ker U$. Of the seventy-two
characteristic vectors satisfying Inequality~\eqref{eq:PartBox} the
following six are the only ones which represent
the given $\SpinC$ structure:
$$
\begin{array}{lll}
(0, 2, 2, 3, 3), &      (0, 0, 0, 3, 3),        &       (0, 2, 2, 1, 1), \\
(0, 0, 0, 1, 1), & (0, 2, 2, -1, -1),   &       (0, 0, 0, -1, -1).
\end{array}$$
Indeed, we claim that of these six characteristic vectors,
$K=(0,0,0,1,1)$ and $-K$ are the only two which can be connected to
characteristic vectors satisfying Inequality~\eqref{eq:OtherPartBox}.
For example, $$(0,2,2,3,3)\sim (4,2,2,-3,-3)\sim U\otimes
(0,4,4,-1,-1).$$ Moreover, since $$\frac{K^2+5}{4}=\frac{3}{4},$$ it
follows immediately that the kernel of $U$ inside $\Combp(G)$ has rank
two, and it is supported in degree $-3/4$.

Next, we claim that the kernel of $U^2$ has rank three. This follows
from the fact that $$U\otimes K \sim U\otimes -K;$$ more specifically:
$$(4,-2,-2,-3,-3)\sim (2,-2,-2,-3,3)\sim (-2,0,0,-1,5)\sim U\otimes
(0,0,0,-1,-1),$$ while $$(4,-2,-2,-3,-3)\sim (-2,2,2,-3,3)\sim
(2,0,0,-5,1)\sim U\otimes (0,0,0,1,1).$$
The restatement for $Y_{12}$ (rather than $-Y_{12}$) follows
from the general properties of the invariant under orientation reversal.
\end{proof}

\begin{prop}
\label{prop:DoubleTrefoil}
Let $Y_0$ denote zero-surgery on the connected sum $T_2$
of two copies of the right-handed trefoil. Then, 
under the identification $\SpinC(Y_0)\cong 2\Z$ (using the first 
Chern class and then a trivialization $H^2(Y_0;\Z)\cong \Z$), we have that
\begin{eqnarray*}
 \HFp(Y_0,\pm 2)&\cong& \Z, \\
 \HFp(Y_0,m)&=& 0 
\end{eqnarray*}
for $|m|>2$. Moreover, as a $\Z[U]$ module, we have that
$$\HFp(Y_0,0)\cong {\mathcal T}^+_{-1/2}\oplus {\mathcal T}^+_{-3/2} \oplus 
\Z_{(-5/2)}$$
(the subscript on the last factor here denotes the absolute grading of the $\Z$ summand).
\end{prop}

\begin{proof}
The fact that $\HFp(Y_0,m)=0$ for $|m|>2$ follows from the adjunction
inequality for $\HFp$ (cf.\ Theorem~\ref{HolDiskTwo:thm:Adjunction}
of~\cite{HolDiskTwo}).  The fact that $\HFp(Y_0,\pm 2)\cong \Z$
follows from the fact that $Y_0$ is a genus two fibered knot
(cf.\ Theorem~\ref{HolDiskSymp:thm:ThreeManifoldsFiber}
of~\cite{HolDiskSymp}).  (An alternative verification of these facts
could be given by a more extensive calculation of $Y_{12}$,
in the spirit of Lemma~\ref{lemma:FirstY12Calc}.)

We now use the graded version of the integral surgeries long exact sequence
(cf.\ Section~\ref{AbsGraded:sec:Lens} of~\cite{AbsGraded}) to determine
$\HFp(Y_0,0)$. Recall that that sequence gives:
$$
\begin{CD}
\cdots@>>>\HFp(S^3)@>{F_1}>> \HFp(Y_0,0) @>{F_2}>> \HFp(Y_{12},Q(0)) @>{F_3}>>\cdots
\end{CD}
$$ (in general, the term involving $Y_0$  reads
$$\bigoplus_{k\in\Z}\HFp(Y_0,12k),$$
but it follows from what we have already seen that $\HFp(Y_0,0)$
is the only non-trivial summand here).  The
map $F_3$ annihilates $\HFinf(Y_{12},Q(0))$, and it can be written as
a sum of terms which decrease degree by at least $-11/4$. Since
$\HFpRed(Y_{12},Q(0))$ is supported in degree $-\frac{1}{4}$, it follows
immediately that $F_3$ is trivial. Bearing in mind that
$F_2$ is a homogeneous map
which shifts degrees by $9/4$
(cf.\ Lemma~\ref{AbsGraded:lemma:CalcDegrees} of~\cite{AbsGraded}),
the result now follows easily.
\end{proof}

\begin{remark}
\label{rmk:OddGens}
It is an easy consequence of this calculation that, if $Y_{-1}$ denotes
the three-manifold obtained by $(-1)$ surgery on the double-trefoil
$T_2$, then $\HFpRed(Y_{-1})$ has generators with both parities;
therefore, so does $\HFpRed(-Y_{-1})$. Note that the plumbing
diagram for $Y_{-1}$ in Figure~\ref{fig:Tref2}
has two bad points (two vertices with degree
three and multiplicity $-1$). This underscores the importance of the
hypothesis on the graph for Theorem~\ref{intro:SomePlumbings}.
\end{remark}

\begin{figure}[ht!]
\relabelbox\small
\cl{\epsfxsize 3in\epsfbox{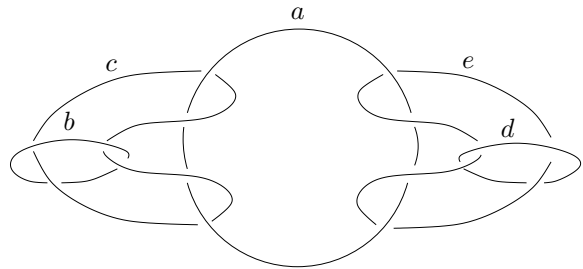}}
\relabel {a}{$a$}
\relabel {b}{$b$}
\relabel {c}{$c$}
\relabel {d}{$d$}
\relabel {e}{$e$}
\endrelabelbox
\caption{\label{fig:GenBorromean}
{\bf{Generalized Borromean rings}}\qua This link has the property that if 
the surgery coefficients $a=b=c=d=e=0$, then the three-manifold obtained is
$S^1\times \Sigma_2$.}
\end{figure}

We now calculate $\HFp(S^1\times \Sigma_2)$. To do this, we think of
$S^1\times \Sigma_2$ as a surgery on a generalized Borromean rings
(compare this with the corresponding calculation of $\HFp(T^3)$ from
Section~\ref{AbsGraded:sec:SampleCalculations}
of~\cite{AbsGraded}). Specifically, consider the link pictured in
Figure~\ref{fig:GenBorromean}. For integers $a,b,c,d,e$, we let
$M(a,(b,c)(d,e))$ denote the three-manifold obtained by surgery
instructions as labelled in the figure (i.e.\ $a$ is the coefficient on
the long circle). In particular, it is easy to see that
$M(0,(1,1)(1,1))$ is zero-surgery on the connected sum of two
right-handed trefoils; while $M(0,(0,0)(0,0))\cong S^1\times
\Sigma_2$.

We calculate $M(0,(0,0),(0,0))$ by successive applications of the long
exact sequence. In this calculation, we will make heavy use of what is
known about $\HFinf$ (see Section~\ref{HolDiskTwo:sec:HFinfty}
of~\cite{HolDiskTwo}). Recall that a three-manifold $Y$ is said to
have standard $\HFinfty$ if for each torsion $\SpinC$ structure
$\spinct_0$, $$\HFinf(Y,\spinct_0)\cong
\Z[U,U^{-1}]\otimes_\Z \Wedge^* H^1(Y;\Z)$$
as a $\Z[U,U^{-1}]\otimes_{\Z}\Wedge^* H_1(Y;\Z)$-module. In general,
we have a spectral sequence whose $E_2$ term is
$\Z[U,U^{-1}]\otimes_\Z \Wedge^* H^1(Y;\Z)$ which converges to
$\HFinf(Y,\spinct_0)$, so this condition is equivalent to the
condition that all higher differentials $d_r$ for $r\geq 2$ are
trivial. Three-manifolds with $b_1(Y)\leq 2$ all have standard
$\HFinf$. Moreover, if $Y$ has standard $\HFinf$, and $\Knot\subset Y$
is a framed, null-homologous knot, and $b_1(Y_{\Knot})=b_1(Y)$, then
$Y_\Knot$ also has standard $\HFinf$
(cf.\ Proposition~\ref{AbsGraded:prop:NegSurgery} of~\cite{AbsGraded}).

Thus, in our exact sequences, we will find it convenient to work with
the three-manifolds $Y=M(a(b,c)(d,e))$ with $a=1$ as much as possible,
since all of these three-manifolds have standard $\HFinf$. The cost is
that these three-manifolds have two extra generators in $\HFp$.
(We will not need to calculate their absolute gradings, however.)

More precisely, we have the following:

\begin{lemma}
\label{lemma:RkTwo}
Fix any integers $b,c,d,e\in \{0,1\}$, and let $Z_0$ denote the three-manifold
$Z_0=M(0,(b,c),(d,e))$.
We have an identification
$$\bigoplus_{\{\spinct\in\SpinC(Z_0)\big| c_1(\spinct)\neq 0\}}\HFp(Z_0,\spinct)\cong \Z^2.$$
Moreover, if $Z_1=M(1,(b,c),(d,e))$, this subgroup injects into $\HFp(Z_1)$.
Similarly, the corresponding subgroup of $\HFa(Z_0)$ has rank four,
and it, too, injects into $\HFa(Z_1)$.
\end{lemma}

\begin{proof}
The first claim follows from the fact that $Z_0$ is a genus two
fibration (see~\cite{HolDiskSymp}).  The injectivity claim follows
from the long exact sequence connecting $Z_0$, $Z_1$, and a third term
which is a connected sum of some number of copies of $S^1\times S^2$.
\end{proof}

We will let $V_{(b,c)(d,e)}\subset \HFp(M(1(b,c)(d,e)))$ denote the
rank two subgroup constructed in Lemma~\ref{lemma:RkTwo}, and
${\widehat V}_{(b,c)(d,e)}\subset \HFa(M(1(b,c)(d,e)))$ be the
corresponding rank four subgroup.\footnote{The Abelian groups which we
meet now and in the rest of this paper are free $\Z$-modules. To
verify this, one can run the exact sequence arguments below for
$\Zmod{p}$ where $p$ is 
an arbitrary prime, and observe that the dimensions of each
of the vector spaces in question is independent of the prime $p$, and
hence, by the universal coefficients theorem, these groups are free
$\Z$-modules. Having said this now, we do not call the readers' attention
to it again.}

\begin{lemma}
\label{lemma:FirstCalc}
We have $\Z[U]$-module identifications:
\begin{eqnarray*}
\HFp(M(1(1,1),(1,1)))&\cong&
        {\mathcal T}^+_{-2}\oplus \Z_{(-3)}\oplus V_{(1,1),(1,1)}, \\
\HFp(M(1(1,0),(1,1)))&\cong&
        {\mathcal T}^+_{-5/2}\oplus {\mathcal T}^+_{-3/2}
\oplus \Z_{(-5/2)}\oplus V_{(1,0),(1,1)}
\end{eqnarray*}
\end{lemma}

\begin{proof} Both are relatively straightforward
applications of the surgery long exact sequence for $\HFp$, given the
calculation from Proposition~\ref{prop:DoubleTrefoil}.

The surgery exact sequence for the triple 
$$
\begin{array}{lll}
\big(M(\infty(1,1),(1,1))\cong S^3, &
M(0(1,1)(1,1)), & M(1(1,1),(1,1))\big)
\end{array}$$ 
reads:
$$
\cdots\to\InjMod{0}\buildrel{F_1}\over\longrightarrow
 \InjMod{-1/2}\oplus\InjMod{-3/2}\oplus \Z_{(-5/2)}
\oplus W \buildrel{F_2}\over\longrightarrow
\HFp(M(1(1,1),(1,1))) \buildrel{F_2}\over\longrightarrow
\cdots$$
where here $W$ is the rank two module generated by the sum of
$\HFp(Y_0,m)$ with $m\neq 0$. Now, we claim that the map $F_3$ is
trivial.  Clearly the map is written as a sum of maps which induce the
trivial map on $\HFinf$ (this is necessary in order for
$\HFinf(M(0(1,1),(1,1)))$ to have its structure). Thus, $F_3$ factors
through $\HFpRed(M(1(1,1)(1,1)))$.  Indeed, it is also trivial on the
image $V_{(1,1)(1,1)}$ of $W$ inside $M(1(1,1)(1,1))$ (by
exactness). But now, since $F_2$ lowers degree by $1/2$, this quotient
group is isomorphic to a single $\Z$ in dimension $-3$. Since the map
$F_3$ does not increase degree, and $\HFp(S^3)$ is supported in
non-negative dimension, it follows that $F_3$ is trivial. Now, $F_1$
is injective, it clearly maps onto the summand $\HFp(Y_0,0)\subset
\HFp(Y_0)$. The first isomorphism claimed in the lemma now follows.

The second isomorphism follows from considering the surgery long exact
sequence for the triple 
$$
\begin{array}{lll}
\big(M(1(1,\infty),(1,1))\cong -\Sigma(2,3,5), 
& M(1(1,0)(1,1)), &
M(1(1,1),(1,1))\big).
\end{array}
$$ Recall that $\HFp(-\Sigma(2,3,5))\cong
\InjMod{-2}$ (see~\cite{AbsGraded}). 
Also, the map in the exact sequence which takes $\HFp(M(1(1,0)(1,1)))$
to $\HFp(M(1(1,1)(1,1)))$ carries the subgroup $V_{(1,0)(1,1)}$ to
$V_{(1,1)(1,1)}$. This latter observation follows from naturality of
the maps induced by cobordisms, togther with the observation that the
corresponding cobordism connecting $M(0(1,0)(1,1))$ to
$M(0(1,1)(1,1))$ induces an isomorphism on the part of $\HFp$
supported in $\SpinC$ structures with non-trivial first Chern
class. This follows easily from the long exact sequence connecting
these $\HFp$ of these three-manifolds.  Indeed, the point here is that
the cobordism from $M(0(1,0)(1,1))$ to $M(0(1,1)(1,1))$ admits a genus
two Lefschetz fibration, see~\cite{HolDiskSymp}, especially
Lemma~\ref{HolDiskSymp:lemma:MaxSpinCIndepOfFib} of that paper.

As before, the map in the long exact sequence taking
$$\HFp(M(1(1,1),(1,1)))\longrightarrow \HFp(-\Sigma(2,3,5))$$ is trivial on the
image of $\HFinf$, and it is also trivial on $V_{(1,1),(1,1)}$ (by
exactness). The remaining quotient group is a $\Z$ in dimension $-3$,
and since the map under consideration does not increase degree, and
$\HFp(-\Sigma(2,3,5))$ is supported in degrees $\geq -2$, it follows
that the map under consideration is trivial.
\end{proof}

\begin{remark}
In view of the above calculations, we see that
$M(1(1,0)(1,1))$ gives yet
another example of an integral homology $S^2\times S^1$ which is
obtained as integral surgery on a two-component link in $S^3$, but which is not
surgery on any knot (indeed, with the given orientation, this manifold
cannot bound an integral homology $D^2\times S^2$). All this
follows from the fact that 
\begin{eqnarray*}
d_{1/2}(M(1(1,0)(1,1)))&=&-\frac{3}{2}, \\
d_{-1/2}(M(1(1,0)(1,1)))&=&-\frac{5}{2}, 
\end{eqnarray*}
together with Theorem~\ref{AbsGraded:thm:IntFormBOneOne}
of~\cite{AbsGraded}.
\end{remark}

\begin{lemma}
We have the identification:
\begin{eqnarray*}
\HFp(M(1(0,0)(0,0))) &{\cong}&  \left({\mathcal T}^+_{0}\right)^6 \oplus 
\left({\mathcal T}^+_{-1}\right)^8 \oplus
\left({\mathcal T}^+_{-2}\right)^2 \oplus V_{(0,0),(0,0)}
\end{eqnarray*}
Moreover, if we let 
$\Z_{(0)}^6\subset \HFa_0(M(1(0,0)(0,0)))$ be 
a subgroup which maps isomorphically onto 
$$(\Ker U|(\InjMod{0})^6)\subset \HFp_{0}(M(1(0,0)(0,0))),$$
then the map induced by the $H_1$-action
$$H_1(M(1(0,0)(0,0));\Z)\otimes \Z_{(0)}^6\longrightarrow 
\HFa_{-1}(M(1(0,0)(0,0)))$$
has six-dimensional image.
\end{lemma}

\begin{proof}
We find it convenient to work with $\HFa$.
The last isomorphism of the previous lemma shows that
$$\HFa(M(1(1,0)(1,1)))\cong \Z^2_{(-5/2)}\oplus \Z^2_{(-3/2)}\oplus {\widehat V}_{(1,0),(1,1)}.$$
Recall  (see Section~\ref{AbsGraded:sec:SampleCalculations} of~\cite{AbsGraded})
that
\begin{eqnarray*}
\HFa(M(1(\infty,0),(1,1)))&\cong& \Z_{(-3/2)}\oplus \Z_{(-5/2)} \\
\HFa(M(1(0,0),(\infty,1))&\cong& \Z^2_{(0)}\oplus \Z^2_{(-1)}.
\end{eqnarray*}
Observe that each $Y=M(1(b,c)(d,e))$ has standard $\HFinf$.  This
gives rise to a subspace of $\HFa(Y)$ of rank $2^{b_1(Y)}$, which
comes from the intersection of the image of $\HFinf(Y)$ in $\HFp(Y)$
with the kernel of $U$. We denote this space ${\widehat W}$.  (Note
that we are dropping $b$, $c$, $d$, and $e$ from the notation
temporarily, as the argument we give here is independent of them).
Moreover, there is also a four-dimensional subspace in $\HFa(Y)$
according to Lemma~\ref{lemma:RkTwo}, which we denote ${\widehat V}$.
We claim that ${\widehat V}\cap{\widehat W}=0$. This follows since (by
construction) ${\widehat W}$ injects into $\HFp(Y)$ under the natural
map $i\colon \HFa(Y)\longrightarrow \HFp(Y)$, but is in the kernel of the
composite of this map with the natural map $\pi\colon \HFp(Y)\longrightarrow
\HFpRed(Y)$. On the other hand, ${\widehat V}$ consists of one two-dimensional
summand which is in the kernel $i$, and another which injects into
$\HFpRed(Y)$ under the composite $\pi\circ i$. Now since
$b_1(M(1(b,c)(d,e)))$ is given by the number of zero entries among the
$b,c,d,e$, the above observations show that
$\HFa$ of each of these three-manifolds satisfies:
\begin{eqnarray*}
\Rk\left(\HFa(M(1(0,1),(1,1)))\right)&\geq& 6 \\
\Rk\left(\HFa(M(1(0,0),(1,1)))\right)&\geq& 8 \\
\Rk\left(\HFa(M(1(0,0),(0,1)))\right)&\geq& 12
\end{eqnarray*}
There are also two  exact sequences, associated
to triples: {\small$$\begin{array}{lll}
\big(M(1(\infty,0),(1,1))\cong -\Sigma(2,3,5)\#(S^2\times S^1), &
M(1(0,0)(1,1)), & M(1(1,1),(1,1))\big) \\
\big(M(1(0,0),(\infty,1)), &
M(1(0,0)(0,1)), & M(1(0,0),(0,0))\big).
\end{array}$$}%
We claim that there are only two possible answers for $\HFa(M(1(0,0)(0,0)))$
which are consistent with all of these constraints:
\begin{equation}
\label{eq:IdentifyHFa}
\HFa(M(1(0,0)(0,0)))\cong \Z_{(0)}^6\oplus \Z_{(-1)}^8\oplus 
\Z_{(-2)}^2\oplus {\widehat V}_{(0,0)(0,0)}.
\end{equation}
or $$\HFa(M(1(0,0)(0,0)))\cong \Z_{(0)}^6\oplus \Z_{(-1)}^9\oplus
\Z_{(-2)}^3\oplus {\widehat V}_{(0,0)(0,0)}.$$
Again, as in the proof of Lemma~\ref{lemma:FirstCalc}, we are using
here the fact that the various ${\widehat V}_{(b,c)(d,e)}$ are mapped
to one another by the corresponding maps, which follows from
naturality of the cobordism invariants, together with the fact that
the relevant cobordisms connecting the corresponding $M(0(b,c)(d,e))$
all admit genus two Lefschetz fibrations.

The latter case is ruled out as follows. 
Suppose it is realized. Then, we consider the long exact sequence
for the triple
{\small$$
\begin{array}{lll}
\big(M(\infty(0,0),(0,0))\cong \#^4(S^2\times S^1), &
M(0(0,0)(0,0))\cong S^1\times \Sigma_2, & M(1(0,0),(0,0))\big).
\end{array}
$$}%
In this case, $\HFa_{3/2}(S^1\times \Sigma_2,\spinct_0)\cong \Z$ (it is the image of
the top-dimensional generator of $\HFa(\#^4(S^2\times S^1))$). Now,
since $\HFa_{-3/2}(S^1\times\Sigma_2)$ surjects onto
$$ \Ker\left(\HFa_{-2}(M(1(0,0),(0,0)))\cong \Z^3 
\longrightarrow \HFa_{-2}(\#^4(S^2\times S^1))\cong \Z\right),
$$
it follows that 
$\Rk\HFa_{-3/2}(S^1\times \Sigma_2)\geq 2$. 
But this contradicts the fact that
$$\HFa_{k}(S^1\times \Sigma_2)\cong \HFa_{-k}(S^1\times\Sigma_2),$$
which follows from the fact that $S^1\times \Sigma_2$ admits an orientation-reversing
diffeomorphism.

It follows that we have isomorphism from Equation~\eqref{eq:IdentifyHFa}, which easily
translates to the claimed identification of $\HFp$.

For the claim about the $H_1$ action, we investigate
the above isomorphisms more carefully. Indeed, we break the verification
into pieces, verifying first that 
the image of
\begin{gather*}
H_1(M(1,(0,0),(0,\infty));\Z)\otimes \HFa_{1/2}(M(1,(0,0),(0,\infty)))
\hspace{1in}
\\
\hspace{2.5in}\longrightarrow \HFa_{-1/2}(M(1(0,0),(0,\infty)))
\end{gather*}
has rank two. But this follows readily from the fact that
$$M(1(0,0),(0,\infty))=M(1(0,0))\# (S^2\times S^1),$$
and thus (according to the K\"unneth principle for connected sums,
cf.\ Proposition~\ref{HolDiskTwo:prop:ConnSum} of~\cite{HolDiskTwo})
$$\HFa(M(1(0,0),(0,\infty))\cong \HFa(M(1(0,0)))\otimes_{\Z} H^*(S^1),$$
where the homology class supported in the $S^2\times S^1$ acts as contraction
on $H^*(S^1)$. In particular, action by this homology class surjects
onto the bottom-dimensional homology of $\HFa(M(1(0,0),(0,\infty)))$
(which in this case is supported in dimension $-1/2$).

We claim also that the map
\begin{gather*}
H_1(M(1,(0,0),(0,1));\Z)\otimes \HFa_{-1/2}(M(1,(0,0),(0,1)))
\hspace{1in}
\\
\hspace{2.5in}\longrightarrow \HFa_{-3/2}(M(1,(0,0),(0,1)))
\end{gather*}
has four-dimensional image. Indeed, we claim that chasing through the
above isomorphisms, the natural maps
\begin{eqnarray*}
\HFinf_d(M(1(0,0),(0,1)))&\longrightarrow& \HFp_d(M(1(0,0),(0,1)))\\
\HFa_d(M(1(0,0),(0,1)))&\longrightarrow& \HFp_d(M(1(0,0),(0,1)))
\end{eqnarray*}
are isomorphisms when $d=-1/2$, $-3/2$. Indeed, the above
natural maps respect the $H_1$-actions. Moreover, since $\HFinf$
of $Y=M(1(0,0),(0,1))$ is standard, and $b_1(Y)=3$, we see that
if $H=H_1(Y;\Z)$, and $H^*=H^1(Y;\Z)$
\begin{eqnarray*}
\HFinf_{-1/2}(Y)&\cong& \Wedge^3 H^* \oplus H^* \\
\HFinf_{-3/2}(Y)&\cong& \Wedge^2 H^* \oplus \Wedge^0 H^*,
\end{eqnarray*}
where the $H=H_1$ action is modelled by contraction. The claim is now
immediate.

Finally, the claim of the lemma then follows
easily from a glance at the $\HFa$-long exact sequence for the triple
$$
\begin{array}{lll}
\big(M(1,(0,0),(0,\infty)), &
M(1,(0,0),(0,0)) &
M(1,(0,0),(0,1))\big),
\end{array}
$$
bearing in mind that all maps are equivariant under the
action of the one-dimensional homology, and using the rank calculations
established above.
\end{proof}

\begin{lemma}
\label{lemma:SmallestHFinf}
Let $Y$ be a three-manifold with $b_1(Y)=5$, and $\spinct_0$ be a
$\SpinC$ structure whose first Chern class is torsion. Then, in each degree
$k$, 
$$9\leq \Rk
\HFinf_k(Y,\spinct_0)\leq 16.$$ Indeed, when the lower bound is
realized, the action $$H_1(Y;\Z)\otimes
\HFinf_{\mathrm{odd}}(Y,\spinct_0) \longrightarrow
\HFinf_{\mathrm{ev}}(Y,\spinct_0)$$ is trivial. (Here, we are using
the absolute $\Zmod{2}$ grading on $\HFinf$, which is characterized by
the property that $\uHFinf_{\ev}$ is non-trivial.)
\end{lemma}

\begin{proof}
Now, to estimate $\HFinf_k(Y,\spinct_0)$, we use the universal coefficients spectral sequence.
The $E_2$ term here is a repeating pattern of
$$
\begin{array}{llllll}
\Z      &       H^1(Y)  &       \Wedge^2 H^1(Y) &       \Wedge^3 H^1(Y) &       \Wedge^4 H^1(Y) &       \Wedge^5 H^1(Y) \\
0       &       0       &       0       &       0       &       0       &       0               
\end{array}
$$ (i.e.\ the repeating pattern comes about by the various $U$ powers).
It is easy to see that the total rank of $E_\infty$ is minimized if
the $d_3$ differential restricts as a surjection from $\Wedge^3
H^1(Y)\longrightarrow \Z,$ an isomorphism from $\Wedge^4
H^1(Y)\longrightarrow H^1(Y)$, and an injection from 
$\Wedge^5 H^1(Y)\longrightarrow 
\Wedge^2 H^1(Y)$.  In that case, the rank of
$\HF_k(Y,\spinct_0)$ (for each $k$) is $9$.  For such a
three-manifold, $\HFinfty(Y)$
is a quotient of a $\Z[U,U^{-1}]\otimes_\Z \Wedge^* H_1(Y;\Z)$-submodule of
$$\left(\Wedge^3 H^1(Y) \oplus
\Wedge^2 H^1(Y)\right)\otimes_\Z \Z[U,U^{-1}].$$ 
In particular, the $H_1(Y)$-action on elements of
odd parity ($\Wedge^2 H^1(Y)$) is trivial. (Recall that the parity is defined so that
$\Wedge^{b_1} H^1(Y)$ has even parity.)
\end{proof}

\begin{theorem}
\label{thm:SOneSigmaTwo}
Letting $\spinct_0\in\SpinC(S^1\times\Sigma_2)$ be the $\SpinC$ structure
with trivial first Chern class, we have the $\Z[U]$-module identification
$$\HFp(S^1\times \Sigma_2,\spinct_0)\cong 
\left({\mathcal T}^+_{3/2}\right) \oplus 
\left({\mathcal T}^+_{1/2}\right)^9
\oplus  \left({\mathcal T}^+_{-1/2}\right)^9 \oplus 
\left({\mathcal T}^+_{-3/2}\right).$$
In particular,
$$\HFa(S^1\times \Sigma_2,\spinct_0)\cong \Z_{(3/2)}\oplus \Z^9_{(1/2)}\oplus
\Z^9_{(-1/2)}\oplus \Z_{(-3/2)}.$$
Moreover, the only other non-trivial $\SpinC$ structures with
non-trivial $\HFp$ are the ones with $c_1(\spinct)=\pm \PD[S^1]$
(where here $[S^1]$ represents the fiber factor of
$S^1\times\Sigma_2$); for each of those, we have that $\HFp$ is
isomorphic to $\Z$.
\end{theorem}

\begin{proof}
Consider the triple: 
$$
\big(M(\infty(0,0),(0,0))\cong \#^4(S^2\times S^1),
M(0(0,0)(0,0))\cong S^1\times \Sigma_2, M(1(0,0),(0,0))\big)
$$ 
Observe that the map
from $\HFp(\#^4(S^2\times S^1))$ to
$\bigoplus_{\{\spinct\big|c_1(\spinct)\neq 0\}}\HFp(\Sigma_2\times S^1)$
is trivial, and hence, so is the map induced on $\HFa$. It follows that
${\widehat W}$ maps isomorphically onto ${\widehat V}$ in the $\HFa$ long-exact sequence
for the triple which now reads:
{\small$$
\cdots
\stackrel{F_1}{\longrightarrow} 
        \HFa(\#^4(S^2\times S^1)) 
\stackrel{F_2}{\longrightarrow} 
        \HFa(S^1\times \Sigma_2,\spinct_0)\oplus{\widehat W} 
\stackrel{F_3}{\longrightarrow} 
        \Z_{(0)}^6\oplus \Z_{(-1)}^8\oplus \Z_{(-2)}^2 \oplus {\widehat V} 
\rightarrow \cdots
$$}%
Of course, as a graded group, we have that 
$$\HFa(\#^4(S^2\times S^1))\cong \Z_{(2)}\oplus
\Z_{(1)}^4 \oplus \Z_{(0)}^6 \oplus \Z_{(-1)}^4\oplus \Z_{(-2)}.$$
Now, let $C$ denote the rank of the kernel of the map
$$F_3\colon \HFa_0(M(1(0,0),(0,0)))\cong \Z^6 \longrightarrow 
\HFa_0(\#^4(S^2\times S^1))\cong \Z^6.$$
Since $\HFa_{k}(S^1\times\Sigma_2)\cong \HFa_{-k}(S^1\times\Sigma_2)$
(which in turn follows from the fact that the three-manifold has an orientation-reversing
diffeomorphism), it follows that 
\begin{equation}
\label{eq:HFaSoneSigTwo}
\HFa(S^1\times \Sigma_2,\spinct_0)=\Z_{(3/2)} \oplus \Z^{4+C}_{(1/2)}\oplus
\Z^{4+C}_{(-1/2)}\oplus \Z_{(-3/2)}.
\end{equation}
Now, it follows from Lemma~\ref{lemma:SmallestHFinf} that $C\geq 4$. 

The case where $C=6$ is excluded by the $H_1$ action as follows.  
We have seen that the map
$$F_2\colon \HFa_{-1}(M(1(0,0),(0,0)))\longrightarrow \HFa_{-1}(\#^4(S^1\times S^2))$$
is surjective; it follows that there must be some element in the six-dimensional
subspace of $\HFa_{-1}(M(1(0,0),(0,0)))$ which is the $H_1$-image of $\HFa_0$
with non-zero projection under $F_2$.
But by the naturality of the $H_1$ action,
such an element must have zero image under $F_2$, since $C=6$ is the hypothesis that
$$F_2\colon \HFa_{0}(M(1(0,0),(0,0)))\longrightarrow \HFa_{0}(\#^4(S^1\times S^2))$$
is identically zero.

To rule out the case where $C=4$, we proceed as follows. If $C=4$,
then the lower bound on the rank of $\HFinf$
Lemma~\ref{lemma:SmallestHFinf} is realized. On the other hand, the
image of the top-dimensional class in $\HFa(\#^4(S^1\times S^2))$ has
odd parity in $\HFa(S^1\times \Sigma_2)$, and yet it has non-trivial
images under the $H_1$ action, contradicting that lemma.

The only remaining case is $C=5$ in Equation~\eqref{eq:HFaSoneSigTwo}. This easily
translates to the claimed identification of $\HFp$.
\end{proof}

\subsection{Further speculation}

Although $\HFp$ of a three-manifold is a subtle invariant, we know that $\HFinf$
is not: it remains unchanged under integral surgeries which preserve $b_1$. 
Still, it is useful to know $\HFinf$ as a starting point for calculations of
$\HFp$. 

As a computational tool, we have a spectral sequence whose
$E_2$ term is given by
$$\Z[U,U^{-1}]\otimes_\Z \Wedge^* H^1(Y;\Z),$$
which converges to $\HFinf(Y)$. Thus a three manifold $Y$ has standard $\HFinf$ if
all the differentials $d_i$ for $i\geq 2$ are trivial.

We have seen (cf.\ Proposition~\ref{AbsGraded:prop:T3}
of~\cite{AbsGraded}) that $T^3$ is a three-manifold whose $\HFinf$ is
not standard. In fact, Theorem~\ref{thm:SOneSigmaTwo} provides us with
another such three-manifold: in each dimension,
the rank of $\HFinf$ (as
a $\Z$-module) is only ten, rather than sixteen.  It is natural to
expect that the cohomology ring of $Y$ plays an important role
here. More concretely, we make the following conjecture (which
is easily seen to be consistent with the above calculations):

\begin{conj}
Let $Y$ be a closed, oriented three-manifold equipped with a torsion $\SpinC$ structure
$\spinct_0$.
The spectral sequence for $\HFinf(Y,\spinct_0)$ collapses after the $E_3$ stage, and moreover
the differential
$$d_3\colon \Wedge^{i}H^1(Y;\Z)\otimes_\Z U^j \longrightarrow \Wedge^{i-3} H^1(Y;\Z)\otimes_\Z
U^{j-1}$$
is given by the homological pairing:
$$d_3(\phi_1\wedge\ldots\wedge\phi_i)
=\frac{1}{3!\cdot(i-3)!}\sum_{\sigma\in{\mathfrak S}_i} (-1)^\sigma \langle 
\phi_{\sigma(1)}\cup\phi_{\sigma(2)}\cup
\phi_{\sigma(3)},[Y]\rangle\cm \phi_{\sigma(4)}\wedge\ldots\wedge\phi_{\sigma(i)},
$$
where ${\mathfrak S}_i$ denotes the permutation group on $i$ letters, and $(-1)^\sigma$
denotes the sign of the permutation.
\end{conj}

\end{document}